\documentclass{article}
\date{}
\usepackage{a4wide}
\usepackage{graphicx}
\newlength\figureheight
\newlength\figurewidth
\usepackage{subfigure}

\usepackage{amssymb}
\usepackage{amsmath}
\usepackage{mathrsfs}
\usepackage{xcolor}
\usepackage{colortbl}
\usepackage{color}

\usepackage{hyperref}
\hypersetup{
    colorlinks,%
    citecolor=black,%
    filecolor=black,%
    linkcolor=black,%
    urlcolor=black
}

\usepackage{mycommands}
\usepackage{multirow}
\usepackage{booktabs}

\usepackage{algorithm}
\usepackage{algorithmic}
\renewcommand{\algorithmicrequire}{\textbf{Input:}}
\renewcommand{\algorithmicensure}{\textbf{Output:}}

\newcommand{\w}{\mathbf{w}}
\newcommand{\z}{\mathbf{z}}

\renewcommand{\xib}{\xi}
\renewcommand{\Bcb}{\Bc}

 \title{A least-squares method for sparse low rank approximation of multivariate functions\thanks{This work was
supported by Airbus Group
and the French National Research Agency (Grant TYCHE
ANR-2010-BLAN-0904).}}

\author{M. Chevreuil\footnotemark[2]\
\and R. Lebrun\footnotemark[3]\ 
\and A. Nouy\footnotemark[2]\ \footnotemark[4]
\and P. Rai\footnotemark[2]}

\begin{document}

\maketitle

\renewcommand{\thefootnote}{\fnsymbol{footnote}}
\footnotetext[2]{Ecole Centrale Nantes, Universit\'e de Nantes,
GeM UMR CNRS 6183, Nantes, France}
\footnotetext[3]{
Department of Applied Mathematics and Modeling, Airbus Group, Suresnes, France
}
\footnotetext[4]{Corresponding author (anthony.nouy@ec-nantes.fr).}

\renewcommand{\thefootnote}{\arabic{footnote}}

\begin{abstract}
In this paper, we propose a low-rank approximation method based on discrete least-squares for the approximation of a multivariate function from random, noisy-free observations. 
Sparsity inducing regularization techniques are used within classical algorithms for low-rank approximation in order to exploit the possible sparsity of low-rank approximations. Sparse low-rank approximations are constructed with a robust updated greedy algorithm which includes an optimal selection of regularization parameters and approximation ranks using cross validation techniques. Numerical examples demonstrate the capability of approximating functions of many variables even when very few function evaluations are available, thus proving the interest of the proposed  algorithm for the 
propagation of uncertainties  through complex computational models.   
  \end{abstract}

\pagestyle{myheadings}
\thispagestyle{plain}
\markboth{Mathilde Chevreuil, R\'egis Lebrun, Anthony Nouy, Prashant Rai}{Least-squares method for sparse low rank approximation}

\section{Introduction}

Uncertainty quantification has emerged as a crucial field of
investigation for various branches of science and engineering. Over
the last decade, considerable efforts have been made in the
development of new methodologies based on a functional point of view
in probability, where random outputs of simulation codes are
approximated with suitable functional expansions.
Typically, when considering a function $u(\xib)$ of input random
parameters $\xib=(\xi_1\ldots \xi_d)$, an approximation is searched
under the form $ u(\xib) \approx \sum_{i=1}^P u_i \phi_i(\xib) $
where the $\phi_i(\xib)$ constitute a suitable basis of
multiparametric
functions (e.g. polynomial chaos basis). 

Several methods have been proposed for the evaluation of functional
expansions (see  \cite{MAT08,NOU09d,LEM10}).
Non intrusive discrete projection methods  allow the
estimation of expansion coefficients by using evaluations
of the numerical model at certain sample points, thus allowing the simple use of existing
deterministic simulation codes.
However the dimension $P$ of classical approximation spaces has an
exponential (or factorial) increase with dimension $d$ and hence the
computational cost becomes prohibitively high as one needs to
evaluate the model for a large number of samples of the order of $P$. The
objective is to construct an approximation of the high
dimensional function $u$, given the fact that we have only limited
information on it. We are particularly interested in the case where
the dimension $d$ is large but the ``effective dimensionality'' of
the function is fairly small.

In order to handle high-dimensional models, we here propose a
low rank tensor approximation method using a discrete least-squares approach, which exploits the
tensor structure of the stochastic function spaces and the possible sparsity of low rank approximations. The underlying
assumption is that the model output functional can be well
approximated using sparse low-rank tensor approximations. 

Low rank approximation methods have
recently been applied to many areas of scientific computing for
approximating elements in high dimensional tensor spaces
\cite{KOL09,HAC12,GRA13,Falco:2012fk,KHO12}, with several applications  in uncertainty propagation \cite{NOU07,DOO09,NOU10,KHO10,MAT11b}. In the context of uncertainty quantification, for
problems involving very high stochastic dimension $d$, instead of
evaluating the coefficients of an expansion in a given approximation
basis (e.g. polynomial chaos), function $u$ is approximated in
suitable low-rank tensor subsets of the form $$
\Mc = \left\{v = F_{\Mc}(\p^{(1)},\ldots,\p^{(d)});\p^{(k)}\in\Real^{n_k}\right\},
$$ where $F_{\Mc}$ is a multilinear map which constitutes a parametrization of the subset $\Mc$
and $\p^{(k)}$ are the parameters. Low rank tensor subsets have nice approximation properties in the sense that they are able to approximate with a good precision a large class of functions that can be encountered in practical applications. The dimension of the parametrization  $\sum_{k=1}^d n_k$ typically grows linearly with
the dimension $d$, thus making possible approximation in high dimension.   Here, we rely on least-squares methods in order to construct approximations in these tensor subsets, using only sample evaluations of the function $u$. 
Methods based on least-squares have already been proposed in  \cite{BEL11,RAI12,DOO12}
 for the construction of tensor
approximations of multivariate functions.  
Here, we propose an
alternative construction that incorporates sparsity-inducing regularization 
allowing the construction
of sparse low rank approximations with only a few function evaluations.
\par
The proposed method consists in approximating the model with a
$m$-term representation $u_m (\xib)= \sum_{i=1}^m \alpha_i w_i
(\xib)$ where the
$\alpha_i$ are real coefficients and where the $w_i$ are successively selected in a sparse low-rank (typically rank-one) tensor subset, ideally $$\Mc^{\boldsymbol{m}\text{-sparse}} = \left\{v = F_{\Mc}(\p^{(1)},\ldots,\p^{(d)});\p^{(k)}\in\Real^{n_k},\Vert \p^{(k)} \Vert_0 \leq m_k\right\},
$$ where $\Vert \cdot\Vert_{{0}}$ is the ``$\ell_0$-norm'' counting the number of non zero coefficients. Although $\Mc^{\boldsymbol{m}\text{-sparse}}$ may have a very small effective dimension 
 $\sum_{k=1}^d m_k\ll \sum_{k=1}^d n_k$, the structure of this set makes optimization in $\Mc^{\boldsymbol{m}\text{-sparse}}$ a combinatorial problem. Therefore, we replace 
 the ideal sparse tensor subset  by 
$$
\Mc^{\boldsymbol{\gamma}} = \left\{v = F_{\Mc}(\p^{(1)},\ldots,\p^{(d)});\p^{(k)}\in\Real^{n_k},\Vert \p^{(k)} \Vert_1 \leq \gamma_k\right\},
$$
where we introduce a convex regularization of the constraints using $\ell_{1}$-norm. 
In practice, optimal approximations in subset $\Mc^{\boldsymbol{\gamma}}$ are computed using 
an alternating minimization algorithm that exploits the specific low dimensional parametrization 
of the subset $\Mc^{\boldsymbol{\gamma}}$ and that involves the solution of successive least-squares problems with sparse $\ell_{1}$-regularization. Cross validation techniques are introduced in order to select optimal regularization parameters   $\boldsymbol{\gamma}$. The progressive (greedy) construction of the $m$-term representation has the advantages of being adaptive and also of reducing the dimension of successive least-squares problems, thus improving the robustness of the least-squares method when only few samples are available.  
As a result, the proposed technique
allows to approximate the response of models with a large number of
random inputs even with a limited number of model evaluations. 
A sparse regularization
technique is also used in order to retain only the most significant functions $w_i$, which results in an improvement of robustness of
the greedy construction when only a 
limited number of samples are available. {\color{black}{The results of this paper highlight the interest of  exploiting both low-rank and sparse structures of functions for a better use of the available information on the function. In this paper, we restrict the presentation to the case where successive corrections are computed in the set $\Mc$ of rank-one tensors. It is well known in practice that greedy rank-one constructions yield suboptimal low-rank canonical decompositions.}}  The extension of the methodology to other low-rank tensor subsets is straightforward. 
 
The outline of the paper is as follows. In section
\ref{sec:functional_approaches}, we introduce some basic concepts
about functional approaches in uncertainty propagation. We also
detail methods based on least-squares for the computation of
approximate functional expansions. In section \ref{sec:tensorreg},
we introduce the proposed sparse low rank approximation method based on
regularized least-squares. Finally the ability of the proposed method
to detect and exploit low rank and sparsity of functions is
illustrated on numerical applications in section \ref{applex}.

\section{Functional representation and least squares methods}
\label{sec:functional_approaches}

\subsection{Stochastic function spaces and their tensor structure}

We here introduce the definitions of stochastic functions spaces and approximation spaces. 
We consider a set $\xib$ of $d$ random variables and we denote by  $(\Xi,\Bc,P_{\xi})$ 
the associated probability space, where $\Xi\subset \Rbb^d$ and where $P_\xi$ is 
the probability law of $\xib$. We suppose that $\xib$ can 
be split into $r$ 
mutually independent sets of random variables 
$\{\xi_k\}_{k=1}^r$, i.e. $\xi=\{\xi_1, \ldots, \xi_r\}$,  where $\xi_k$ takes 
values in  $\Xi_k\subset \Rbb^{d_k}$. We have $d=\sum_{k=1}^r d_k$.  We denote by 
$(\Xi_k,\Bc_k,P_{\xi_k})$
the probability space associated with $\xi_k$, where $P_{\xi_k}$ is the probability law of
$\xi_k$. 
Therefore, the probability space $(\Xi,\Bc,P_{\xi})$
associated with $\xib=(\xi_1,\hdots,\xi_r)$ has a
product structure with $ \Xi=\times_{k=1}^r\Xi_k$ and  $P_{\xi}=\otimes_{k=1}^r P_{\xi_k} $.

We denote by $L^2_{P_\xib}({\Xi})$ the Hilbert space of second order
random variables defined on $(\Xi,\Bcb,P_\xib)$, defined by 
$$
L^2_{P_\xib}({\Xi}) = \left\{ u : y\in\Xi \mapsto u(y) \in\Rbb ; \int_{\Xi} u(y)^2 P_{\xi}(dy)<\infty \right\},
$$ 
which is a tensor
Hilbert space with the following tensor structure:
$$
L^2_{P_\xib}({\Xi}) = L^2_{P_{\xi_1}}(\Xi_1) \otimes \hdots \otimes
L^2_{P_{\xi_r}}(\Xi_r).
$$
We now introduce approximation spaces $\Sc^k_{n_k} \subset
L^2_{P_{\xi_k}}(\Xi_k)$ with orthonormal basis $\{\phi_{j}^{(k)}\}_{j=1}^{n_k}$, such that
$$
\Sc_{n_k}^k = 
\left\{v^{(k)}(y_k) = \sum_{j=1}^{n_k} v_j^{k}
\phi_{j}^{(k)}(y_k);v_j^k\in\Rbb\right\} 
= \left\{v^{(k)}(y_k) 
= \phib^{(k)}(y_{k})^{T} \v^{(k)}
;\v^{(k)}\in\Rbb^{n_{k}}\right\} ,
$$
where $\v^{(k)}$ denotes the vector of coefficients of $v^{(k)}$ and 
where $\phib^{(k)}=(\phi_1^{(k)},\ldots ,\phi_{n_k}^{(k)})^T$ denotes the vector of basis functions.
 An approximation space $\Sc_{\boldsymbol{n}}\subset
 L^2_{P_\xib}({\Xi})$ is then obtained by tensorization of approximation spaces $\Sc^k_{n_k} $:
\begin{align*}
\Sc_{\boldsymbol{n}} &= \Sc_{n_1}^1\otimes \hdots \otimes \Sc_{n_r}^r
=\left\{v = \sum_{i \in I_{\boldsymbol n}} 
v_{i} \phi_{i} \; ; \; v_{ i} \in\Rbb\right\},
\end{align*}
where $I_{\boldsymbol n} =  \times_{k=1}^r
\{1 \hdots n_k\}$ and 
$\phi_{ i}(y)=(\phi_{i_1}^{(1)} \otimes \hdots \otimes
\phi_{i_r}^{(r)})(y_1,\hdots,y_r) = \phi_{i_1}^{(1)}(y_1)
\hdots  \phi_{i_r}^{(r)}(y_r)$. An element
$v = \sum_{i} v_{i} \phi_{i} \in\Sc_{\boldsymbol{n}}$ can be identified with the algebraic
tensor $\mathbf{v}\in  \Rbb^{n_1}\otimes
\hdots\otimes \Rbb^{n_r}$ such that $ (\v)_{ i} = v_{ i}$. 
Denoting $\phib(y)=\phib^{(1)}(y_1)\otimes\ldots\otimes\phib^{(r)}(y_r) \in \Rbb^{n_1}\otimes
\hdots\otimes \Rbb^{n_r} $, we have the identification $\Sc_{\boldsymbol{n}} \simeq \Rbb^{n_1}\otimes
\hdots\otimes \Rbb^{n_r}$ with 
$$
\Sc_{\boldsymbol{n}} = \left\{ v(y) = \langle \boldsymbol{\phi}(y),\v \rangle ; \v\in \Rbb^{n_1}\otimes
\hdots\otimes \Rbb^{n_r}\right\},
$$
where $\langle\cdot,\cdot\rangle$ denotes the canonical inner product in $\Rbb^{n_1}\otimes
\hdots\otimes \Rbb^{n_r}$.

Here, we suppose that the approximation space $\Sc_{\boldsymbol n}$ is
given and sufficiently rich to allow accurate representations of a
large class of functions (e.g. by choosing polynomial spaces with
high degree, wavelets with high resolution...). Then, the aim is to provide a method for the approximation
of functions in $\Sc_{\boldsymbol n}$ for high dimensional
applications and using only limited information on the functions. 
Note that in the case $r=d$, approximation space
$\Sc_{\boldsymbol{n}}$ has a dimension $ \prod_{k=1}^d n_k$ which
grows exponentially with the dimension $d$, thus making impossible
the numerical representation and computation of an element
$v\in\Sc_{\boldsymbol{n}}$ for high dimensional applications.
In order to reduce the dimensionality, a classical strategy consists in introducing approximation subspaces $\Sc_{\boldsymbol{n},p}\subset \Sc_{\boldsymbol{n}}$ that are constructed using suitable tensorization rules: $\Sc_{\boldsymbol{n},p} = 
\{v = \sum_{i \in I_{\boldsymbol{n},p}} 
v_{i} \phi_{i} \; ; \; v_{ i} \in\Rbb\},$
where $I_{\boldsymbol{n},p} \subset I_{\boldsymbol n}$ is an index set which can be chosen a priori. For $r=d$, a
typical construction consists in taking for $\Sc_{n_k}^{k}$ the
space of degree $p$ polynomials $\Pbb_p(\Xi_k)$, with $\phi^{(k)}_j$ the orthogonal polynomial of degree $j-1$, and for $I_{\boldsymbol{n},p} =
\{{i}\in I_{\boldsymbol n}; \sum_{k=1}^d (i_k-1) \le
p\}$. Thus, $\Sc_{\boldsymbol{n},p} $ appears to be the so called polynomial chaos
composed of multidimensional polynomials with total degree less than
$p$ \cite{GHA91,XIU02}. 
Other tensorization strategies have been proposed, based on a priori knowledge on the solution or based on adaptive strategies. 
In the present work, we are interested in alternative methods that try to approximate high dimensional functions using low-rank approximations. They will be introduced in section \ref{sec:tensorreg}.

\subsection{Least-squares methods}
We here consider the case of a real-valued model output $u:\Xi
\rightarrow \Rbb$. We denote by $ \{y^q\}_{q=1}^Q\subset \Xi$ a set
of $Q$ samples of $\xib$, and by $\{u(y^q)\}_{q=1}^Q\subset \Real$
the corresponding function evaluations. We suppose that an
approximation space $\Sc_P = span\{\phi_i\}_{i=1}^P$ is given.
Classical least-squares method for the construction of an
approximation $u_P\in\Sc_P$ then consists in solving the following
problem:
\begin{align}
\Vert u-u_P\Vert_Q^2 = \min_{v\in\Sc_P} \Vert u-v\Vert_Q^2
\quad\text{with} \quad\Vert u \Vert_Q^2 = \frac{1}{Q} \sum_{q=1}^Q 
u(y^q)^2. \label{eq:least-squares_problem}
\end{align}
Note that $\Vert\cdot\Vert_Q$ only defines a semi-norm on
$L^2_{P_\xi}(\Xi)$ but it may define a norm on the finite
dimensional subspace $\Sc_{P}$ if we have a sufficient number $Q$ of
model evaluations. A necessary condition is $Q\ge P$. However, this
condition may be unreachable in practice for high dimensional
stochastic problems and usual a priori (non adapted) construction of
approximation spaces $\Sc_P$. Moreover, classical least-squares method may
yield bad results because of ill-conditioning (solution very
sensitive to samples). A way to circumvent these issues is
to introduce a regularized least-squares functional:
\begin{align}
\Jc^\lambda(v) = \| u-v\|_Q^2 + 
\lambda\Lc(v),\label{eq:regularized_least-squares_functional}
\end{align}
where $\Lc$ is a
regularization functional and where $\lambda$ refers to some regularization parameter. The regularized least-squares problem then
consists in solving
\begin{align}
\Jc^\lambda(u_P^\lambda) = \min_{v\in\Sc_P}
\Jc^\lambda(v).\label{eq:regularized_least-squares_problem}
\end{align}
Denoting by $\v=(v_1,\hdots,v_P)^T\in\Rbb^P$ the coefficients of an
element $v = \sum_{i=1}^P v_i \phi_i \in \Sc_P$, we can write
\begin{align}
\Vert u-v\Vert_Q^2 = \Vert \mathbf{z}-\Phib \v\Vert_2^2,
\end{align}
with $\mathbf{z} = (u(y^1),\hdots,u(y^Q))^T \in\Rbb^Q$ the vector of
random evaluations of $u(\xib)$ and $\Phib \in \Real^{Q\times P}$
the matrix with components $ (\Phib)_{q,i} = \phi_i(y^q)$. We can
then introduce a function $L:\Rbb^P\rightarrow \Rbb$ such that
$\Lc(\sum_{i}v_i \phi_i) = L(\v)$, and a function
$J^\lambda:\Rbb^P\rightarrow \Rbb$ such that
$\Jc^\lambda(\sum_{i}v_i \phi_i) = J^\lambda(\v)= \Vert \mathbf{z} -
\Phib \v\Vert_2^2 + \lambda L(\v) $.  An algebraic version of
least-squares problem \eqref{eq:regularized_least-squares_problem} can
then be written as follows:
\begin{align}
\min_{\v \in \Rbb^P} \Vert \mathbf{z}-\Phib \v\Vert_2^2 + 
\lambda L(\v).\label{eq:regularized_least-squares_problem_algebraic}
\end{align}
Regularization introduces additional information such as smoothness,
 sparsity, etc. Under some assumptions on the
regularization functional $\Lc$, problem
\eqref{eq:regularized_least-squares_problem} may have a unique
solution. However, the choice of regularization strongly influences
the quality of the obtained approximation. {Another significant 
component of solving \eqref{eq:regularized_least-squares_problem_algebraic} is the choice 
of regularization parameter $\lambda$. In this paper, we use cross validation for the
selection of an optimal value of $\lambda$.}

\subsection{Sparse regularization}\label{sparse_reg}
Over the last decade, sparse approximation methods  have
been extensively studied in different scientific disciplines. 
 A sparse function is one that can be represented using few non zero terms when expanded on a suitable
basis. In the context of uncertainty quantification, if a stochastic
function is known to be sparse on a particular function basis, e.g.
polynomial chaos (or tensor basis), sparse regularization methods
can be used for quasi optimal recovery with only few sample
evaluations. In general, {\color{black}a successful reconstruction of sparse
solution vector depends on  sufficient sparsity of the
coefficient vector and additional properties (incoherence) depending on the samples and of the chosen basis (see \cite{CAN06b,DON06} or \cite{DOO11} in the context of uncertainty quantification).
This strategy has been found to be effective for
non-adapted sparse approximation of the solution of some PDEs \cite{BLA11,DOO11}.
}
\\\par
More precisely, an approximation $ \sum_{i=1}^P u_i\phi_i(\xi)$ of a
function $u(\xi)$ is considered as sparse on a particular basis
$\{\phi_i(\xib)\}_{i=1}^P$ if it admits a good approximation with only a few 
non zero coefficients. Under certain conditions,
a sparse approximation can be computed accurately using only
$Q \ll P$ random samples of $u(\xib)$ via sparse regularization.

Given the random samples $\mathbf{z}\in\Rbb^Q$ of the function $u(\xib)$, a best $m$-sparse (or $m$-term) approximation 
of $u$ can be ideally obtained by
solving the constrained optimization problem 
\begin{align}
\min_{\v\in\Rbb^P}  \|{ \mathbf{z}-\Phib\v }\|_2^2 \quad \mathrm{subject\;to}
\quad \Vert{\v}\Vert_0 \le m , \label{eq:min_P0delta}
\end{align}
 where $\norm{\v}_0 = \#\{i\in
\{1,\ldots,P\}\;:\;v_i\neq 0\}$ is the so called $\ell_{0}$-``norm'' of $\v$ which gives the number of non zero components
of $\v$.
Problem \eqref{eq:min_P0delta} is a combinatorial optimization problem which is NP hard to solve. Under certain assumptions, problem
\eqref{eq:min_P0delta} can be reasonably well approximated by
the  following constrained optimization problem which introduces a convex relaxation of the $\ell_{0}$-``norm'':
\begin{align}
\min_{\v\in\Rbb^P} \|{ \mathbf{z}-\Phib\v }\|_2^2 \quad \mathrm{subject\;to} \quad  \|{\v}\|_1
 \le \delta, 
\label{eq:min_P1delta}
\end{align}
where $\norm{\v}_1 = \sum_{i=1}^P |v_i|$ is the $\ell_1$-norm of $\v$. Since the $\ell_2$ and $\ell_1$-norms are convex, we can equivalently consider the following convex optimization problem, known as Lasso \cite{TIB96} or  basis pursuit \cite{CHE99}:
\begin{align}
\min_{\v \in \Rbb^P} \|{\mathbf{z}-\Phib\v}\|_2^2+\lambda\|\v\|_1,
\label{eq:min_P1lambda}
\end{align}
where $\lambda>0$ corresponds to  a Lagrange multiplier whose value is related to 
$\delta$. Problem \eqref{eq:min_P1lambda} appears as a regularized least-squares problem. The $\ell_{1}$-norm is a sparsity inducing regularization function in the sense that the solution $\v$ of  \eqref{eq:min_P1lambda} may contain components which are exactly zero.  Several optimization algorithms have been
proposed for solving \eqref{eq:min_P1lambda}
(see \cite{BAC12}).  In this paper, we use the Lasso modified least angle 
regression algorithm (see LARS presented in \cite{EFR04}) that provides a set of
$N_r$ solutions, namely the regularization path, with increasing $\ell_1$-norm. 
Let ${\v^j}$, with ${j=1, \ldots, N_r}$, denote this set of solutions, ${A_j}\subset \{1,\hdots,P\}$ be
the index set  corresponding to non zero coefficients of $\v^{j}$, $\v^j_{A_{j}} \in \Rbb^{\#A_{j}}$ the vector of the coefficients $A_j$ of $\v^j$, and $\Phib_{A_j} \in \Rbb^{Q\times \#A_{j}}$ the submatrix of $\Phib$ obtained by extracting the columns of $\Phib$  corresponding to indices $A_{j}$.
The optimal solution $\v$ is then selected 
using the fast leave-one-out cross validation error estimate \cite{CAW04} which relies on the use of the Sherman-Morrison-Woodbury formula (see  \cite{BLA11} for its implementation within Lasso modified LARS algorithm).
Algorithm \ref{alg:loo_cv} briefly outlines the cross validation procedure for the selection of the optimal solution. In this work, we have used Lasso modified LARS implementation 
of SPAMS software \cite{MAI10} for $\ell_1$-regularization.

%

\begin{algorithm}[H]
\caption{Algorithm to determine optimal LARS solution using leave-one-out cross validation.}
\label{alg:loo_cv}
\begin{algorithmic}[1]
\REQUIRE sample vector $\mathbf{z} \in\Rbb^Q$ and matrix $\Phib \in \Real^{Q\times P}$
\ENSURE vector of coefficients $\v\in \Rbb^P$
\STATE Run the Lasso modified LARS procedure to obtain $N_{r}$ solutions ${\v^1,\hdots,\v^{N_{r}}}$ of the regularization path, with corresponding sets of non zeros coefficients $A_1,\hdots,A_{N_r}$.
\FOR {$j=1,\ldots, N_r$}
\STATE Recompute the non zero coefficients $\v^{j}_{A_j}$ of ${\v^{j}}$ using ordinary least-squares: \\ $\v_{A_j}^{j}=\arg\min_{\v \in \Rbb^{\#A_j}} \|{\mathbf{z}-\Phib_{A_j}\v}\|_2^2$
\STATE Compute $h_q=(\Phib_{A_j}(\Phib_{A_j}^T\Phib_{A_j})^{-1}\Phib_{A_j}^T)_{qq}$.
\STATE Compute relative leave-one-out error 
$\epsilon_{j}=\frac{1}{Q}\sum_{q=1}^Q\Big(\frac{(\z)_q-(\Phib_{A_j}\v^j_{A_j})_q}{(1-h_q)\hat\sigma(\z)}\Big)^2$, where 
$\hat\sigma(\z)$ is the empirical standard deviation of $\z$.
\ENDFOR
\STATE Return the optimal solution $\v$ such that $\v_{A_{j^\ast}}=\v^{j^\ast}_{A_{j^\ast}}$ with $j^{\ast}=\arg\min_j{\epsilon_j}$.
\end{algorithmic}
\end{algorithm}

\section{Sparse low-rank tensor approximations based on least squares method}
\label{sec:tensorreg}

The aim is to find a sparse low rank approximation of a function $u(\xi)$ in the finite dimensional tensor space $\Sc_{\boldsymbol{n}} =
\Sc^1_{n_1}\otimes \hdots \otimes \Sc^r_{n_r}$. The proposed method relies on a greedy algorithm where successive corrections of the approximation are computed in small low-rank tensor subsets using least-squares with sparsity inducing regularization.
Here, we restrict the presentation to the case where successive corrections are computed in the elementary set of rank-one tensors $\Rc_1$, thus resulting in the construction of a low-rank canonical approximation of the solution. 
However, the methodology could be naturally extended in order to construct sparse representations in other tensor formats.  

\subsection{Sparse canonical tensor subsets}
 Let $\Rc_1 $ denote the set 
of (elementary) rank-one tensors in $\Sc_{\boldsymbol{n}} =
\Sc^1_{n_1}\otimes \hdots \otimes \Sc^r_{n_r}$, defined by 
\begin{align*}
\Rc_1 &= \left\{w(y) = \left(\otimes_{k=1}^r w^{(k)}\right)(y) =
\prod_{k=1}^r w^{(k)}(y_k) \; ; \; w^{(k)} \in \Sc_{n_k}^k\right\},
\end{align*}
or equivalently by
\begin{align*}
\Rc_1 = \left\{w(y) = \langle\phib(y),{\w^{(1)}}\otimes\ldots \otimes{\w^{(r)}}\rangle;\w^{(k)}\in\Real^{n_k}\right\},
\end{align*}
where $\phib(y) = \phib^{(1)}(y_1)\otimes\hdots\otimes \phib^{(r)}(y_r)$, with $\phib^{(k)}=(\phi_1^{(k)},\ldots ,\phi_{n_k}^{(k)})^T$ the vector of basis functions of  $\Sc^{k}_{n_k}$, and where 
$\w^{(k)}=(w_1^k,\hdots,w_{n_k}^k)^T$ is the set of coefficients of $w^{(k)} $ in the basis of $ \Sc_{n_k}^k$, that means 
$w^{(k)}(y_k) = \sum_{i=1}^{n_k} w_i^k\phi_i^{(k)}(y_k) $.


Approximation in $\Rc_1$ using
classical least-squares methods possibly enables to recover a good 
approximation of the solution using a reduced number of samples. However, 
it may not be sufficient in the case where the approximation 
spaces $\Sc^k_{n_k}$ have high dimensions $n_k$, thus resulting 
in a manifold of rank-one elements $\Rc_1$ with high dimension $\sum_{k=1}^r n_k$.
This difficulty may be circumvented by introducing approximations in a $\boldsymbol{m}$-sparse rank-one subset defined as 
\begin{align*}
\Rc_1^{\boldsymbol{m}\text{-sparse}} = \left\{w(y) = \langle \phib(y),{\w^{(1)}}\otimes\ldots \otimes{\w^{(r)}}\rangle;\w^{(k)}\in\Real^{n_k},\Vert \w^{(k)} \Vert_0 \leq m_k\right\}
\end{align*}
with effective dimension $\sum_{k=1}^r m_k\ll \sum_{k=1}^r n_k$. As mentioned in section \ref{sparse_reg}, performing 
least-squares approximation in this set may not be computationally tractable. We thus introduce a convex relaxation of the $\ell_0$-``norm''
to define the subset $\Rc_1^\gammab$ of $\Rc_1$ defined as
\begin{align*}
\Rc_1^{{\gammab}} =\left\{w(y) = \langle \phib(y),{\w^{(1)}}\otimes\ldots \otimes{\w^{(r)}}\rangle;\w^{(k)}\in\Real^{n_k},\Vert \w^{(k)} \Vert_1 \leq \gamma_k\right\},
\end{align*}
where the set of parameters $(\w^{(1)},\hdots,\w^{(r)})$ is now searched in a convex subset of $\Rbb^{n_1}\times \hdots \times \Rbb^{n_r}$. 
\\
Finally, we introduce the set of canonical rank-$m$ tensors  $\Rc_m = \left\{ v = \sum_{i=1}^m w_i \;; w_i\in\Rc_1  \right\} $ and the corresponding subset
\begin{align*}
\Rc_m^{\boldsymbol{\gamma}^1,\hdots,\boldsymbol{\gamma}^m}  = \left\{ v = \sum_{i=1}^m w_i \;; w_i\in\Rc_1^{\boldsymbol{\gamma}^i} \right\}.
\end{align*}
In the following, we propose algorithms for the
construction of approximations in tensor subsets $\Rc_1^{\boldsymbol{\gamma}}$ and
$\Rc_m^{\boldsymbol{\gamma}^1,\hdots,\boldsymbol{\gamma}^m} $. These subsets are low-dimensional subsets of the approximation
space $\Sc_{\boldsymbol{n}}$ but are not linear spaces nor
convex sets, thus making more difficult the analysis and practical
resolution of optimization problems in these sets. In practice, we rely on heuristic alternating minimization algorithms presented in the following sections.
\begin{remark}\label{rem:rmclosed}
For the construction of a rank-$m$ approximation, we will introduce a progressive construction based on successive rank-one corrections. A direct approximation in $\Rc_m$ would also be possible, with a straightforward extension of the algorithm presented  in Section \ref{sec:updategreedy}. 
Best approximation problem in $\Rc_m$ for $m\ge 2$  is an ill-posed problem since $\Rc_m$ is not closed (see e.g. \cite{HAC12} lemma 9.11 pg 255). However, using sparsity-inducing regularization (and also other types of regularizations)  makes the best approximation problem well-posed. Indeed, it can be proven that the subset $\Rc_m^{\boldsymbol{\gamma}^1,\hdots,\boldsymbol{\gamma}^m} $ of canonical tensors with bounded factors is a closed subset. 
\end{remark}

\begin{remark}
Other tensor subsets have been introduced that have better
approximation properties, such as Tucker tensor sets or Hierarchical Tucker
tensor sets {\color{black}(see \cite{HAC12} for a comprehensive review)}. 
These tensor formats are not
considered here.
\end{remark}

\subsection{Construction of sparse rank-one tensor approximation}\label{sec:updategreedy}

The subset $\Rc_1^{\boldsymbol \gamma}$ can be parametrized with the set
 of parameters $(\w^{(1)},\hdots,\w^{(r)})\in \Rbb^{n_1}\times\hdots\times \Rbb^{n_r}$ such that
$\Vert \w^{(k)}\Vert_1\leq \gamma_k$ ($k=1,\ldots, r$), 
this set of parameters corresponding to an element $w = \otimes_{k=1}^r w^{(k)}$ where  
$\w^{(k)}$ denotes the vector of coefficients of an 
element $w^{(k)}$. 
With appropriate choice of bases, a sparse rank-one function $w$ could be well approximated 
using vectors $\w^{(k)}$ with only a few non zero coefficients. 

We compute a rank-one approximation $w=\otimes_{k=1}^r
w^{(k)}\in \Rc_1^{\boldsymbol\gamma}$ of $v$ by solving 
the least-squares problem 
\begin{align}
\min_{w\in\Rc_1^\gammab} \Vert v-w\Vert_Q^2 = \min_{\substack{\w^{(1)}\in\Rbb^{n_1},\hdots,\w^{(r)}\in\Rbb^{n_r}\\\Vert \w^{(1)}\Vert_1\le \gamma_1,\hdots,\Vert \w^{(r)}\Vert_1\le \gamma_r
}} \Vert v-\langle \phib,\w^{(1)} \otimes \hdots \otimes \w^{(r)}\rangle \Vert_Q^2.\label{eq:sparse_rank1_correction}
\end{align}
Problem \eqref{eq:sparse_rank1_correction} can be equivalently written 
\begin{align}
\min_{ \w^{(1)}\in\Rbb^{n_1},\hdots,\w^{(r)}\in\Rbb^{n_r}} 
\Vert v-\langle \phib,\w^{(1)} \otimes \hdots \otimes \w^{(r)}\rangle \Vert_Q^2 + \sum_{k=1}^r \lambda_k \Vert \w^{(k)}\Vert_1,\label{eq:sparse_rank1_correction_bis}
\end{align}
where the values of the regularization parameters $\lambda_k>0$ (interpreted as Lagrange multipliers) are related to $\gamma_k$. 
In practice, minimization problem \eqref{eq:sparse_rank1_correction_bis}
is solved using an alternating minimization algorithm which consists in successively minimizing over 
$\w^{(j)}$ for  fixed values of $\{\w^{(k)}\}_{k\neq
j}$. Denoting by $\mathbf{z} \in \Rbb^Q$ the vector of samples of function $v(\xi)$ and by $\Phib^{(j)} \in \Rbb^{Q\times n_j}$ the matrix whose
components are $(\Phib^{(j)})_{qi} = \phi_i^j(y^q_j) \prod_{k\neq
 j}w^{(k)}(y_k^q)$, the minimization problem on  $w^{(j)}$ can be written
 \begin{align}
\min_{\w^{(j)}\in\Rbb^{n_j}} \|{\mathbf{z}-\Phib^{(j)}\w^{(j)}}\|_2^2 + \lambda_j\Vert\w^{(j)}\Vert_1, \label{eq:sparse_rank1_correction_alternativemin_regul_discrete}
\end{align}
which has a classical form of a least-squares problem with a sparsity inducing $\ell_1$-regularization. 
Problem \eqref{eq:sparse_rank1_correction_alternativemin_regul_discrete} is solved using the Lasso modified LARS algorithm where the optimal solution is selected using the leave-one-out cross validation procedure presented in Algorithm \ref{alg:loo_cv}. Algorithm \ref{sparse_rank_one_algo} outlines the construction of a sparse rank one approximation.
\begin{algorithm}[h!]
\caption{Algorithm to compute sparse rank one approximation of a function $v$.}
\label{sparse_rank_one_algo}
\begin{algorithmic}[1]
\REQUIRE vector of evaluations $\z = (v(y^1),\hdots,v(y^Q))^T \in\Rbb^Q$.
\ENSURE rank-one approximation $w(y)=\langle \phib(y), \w^{(1)},\hdots,\w^{(r)}\rangle$.
\STATE Initialize the vectors $\{\w^{(k)}\}_{k=1}^r$ and set $l=0$. 
\STATE \label{algo2-loop} $l \leftarrow l+1$.
\FOR {$j=1,\ldots, r$}
\STATE Evaluate matrix $\Phib^{(j)}$.
\STATE\label{algo2-stepals} Solve problem \eqref{eq:sparse_rank1_correction_alternativemin_regul_discrete} using Algorithm \ref{alg:loo_cv} for input $\z\in \Rbb^Q$ and $\Phib^{(j)} \in \Rbb^{Q\times n_j}$ to obtain $\w^{(j)}$. \label{step:alg2_l1}
\ENDFOR
\STATE Compute   $\hat{\z}=(w(y^1),\hdots,w(y^Q))^T $.
\IF{$\Vert \z - \hat{\z} \Vert_2 > \epsilon$ and $l\leq l_{max}$}
\STATE Go to Step \ref{algo2-loop}.
\ENDIF
\STATE Return the solution parameters $\w^{(1)},\hdots,\w^{(r)}$.
\end{algorithmic}
\end{algorithm}

\begin{remark}[Other types of regularization]\label{algo_versions}
Different rank-one approximations can be defined by changing the type of regularization and constructed  by replacing 
step \ref{algo2-stepals} of Algorithm  \ref{sparse_rank_one_algo}. 
First, one can consider 
a simple least-squares without regularization, namely Ordinary Least Squares (OLS), by replacing 
 step \ref{algo2-stepals}  by the solution of
\begin{align}
\min_{\w^{(j)}\in\Rbb^{n_j}} \|{\mathbf{z}-\Phib^{(j)}\w^{(j)}}\|_2^2 .\label{eq:sparse_rank1_correction_alternativemin_noregul_discrete}
\end{align}
Also, one can consider a regularization using $\ell_2$-norm (ridge regression) by replacing step \ref{algo2-stepals} by the solution of 
\begin{align}
\min_{\w^{(j)}\in\Rbb^{n_j}} \|{\mathbf{z}-\Phib^{(j)}\w^{(j)}}\|_2^2+  \lambda_j\Vert\w^{(j)}\Vert_2^2\label{eq:sparse_rank1_correction_alternativemin_regul2_discrete}
\end{align}
with a selection of optimal parameter $\lambda_j$ using standard cross-validation (typically  $k$-fold cross-validation).
The approximations obtained with these different variants will be  compared in the numerical examples of section \ref{applex}.
\end{remark}

\subsection{Updated greedy construction of sparse rank-$M$ approximation}
\label{updated_greedy_const}
We now wish to construct a sparse rank-$M$ approximation $u_M\in\Rc_M$ of $u$ of the form
$ u_M = \sum_{m=1}^M \alpha_m w_m \label{eq:rank_m_approximation}$ by successive computations of sparse rank-one approximations $w_m=\otimes_{k=1}^r w_m^{(k)}$.  
{\begin{remark}This construction yields a suboptimal rank-$M$ approximation but it has several advantages: successive minimization problems in $\Rc_1$ are well-posed (without any regularization), it requires the learning of a small number of parameters at each iteration. However, let us recall that a direct approximation in sparse canonical tensor format would also be possible (see remark \ref{rem:rmclosed}).\end{remark}}

We start by setting $u_0=0$. Then, knowing an approximation $u_{m-1}$ of $u$, we proceed as follows. 
\subsubsection{Sparse rank-$1$ correction step}
We first compute a
correction $w_m\in \Rc_1 $ of $u_{m-1}$ by solving 
\begin{align}
\min_{w\in\Rc_1^{\gammab}} \Vert u-u_{m-1}-w\Vert_Q^2, 
\label{eq:sparse_greedy_correction}
\end{align}
which can be reformulated as 
\begin{align}
\min_{w\in\Rc_1 } \Vert u-u_{m-1}-\langle \phib,\w^{(1)}\otimes \hdots \otimes \w^{(r)}\rangle \Vert_Q^2 +  \sum_{k=1}^r \lambda_k \Vert\w^{(k)}\Vert_1.
\label{eq:sparse_greedy_correction_regul}
\end{align}
Problem \eqref{eq:sparse_greedy_correction_regul}
is solved using an alternating minimization algorithm, which consists in successive minimization problems of the form 
 \eqref{eq:sparse_rank1_correction_alternativemin_regul_discrete}  where $\mathbf{z} \in \Rbb^Q$ is the vector of samples of the 
residual $(u-u_{m-1})(\xib)$.  Optimal parameters $\{\lambda_k\}_{k=1}^r$ are selected with the fast leave-one-out cross-validation. 

\subsubsection{Updating step}

Once a rank-one correction $w_m$ has been computed, it is normalized and the approximation
 $u_m = \sum_{i=1}^m \alpha_i w_i$ is computed 
by solving a regularized least-squares problem:
\begin{align}
\min_{\alphab=(\alpha_1,\hdots,\alpha_m) \in \Rbb^m} \Vert u-\sum_{i=1}^m\alpha_i w_i\Vert_Q^2 + 
\lambda' \Vert \alphab \Vert_1 . \label{eq:algo_update_regularized}
\end{align}
This updating step allows a selection of
significant terms in the canonical decomposition, that means when some $\alpha_i$
are found to be negligible, it yields an approximation 
$u_m = \sum_{i=1}^m \alpha_i w_i $ with a lower effective
rank representation. The selection of parameter $\lambda'$ is also done with a cross-validation technique. 
\begin{remark}\label{rem:tucker_update}
Note that an improved updating strategy could be introduced as follows. 
At step $m$, denoting by $w_i = \otimes_{k=1}^r w_i^{(k)}$, $1\le i \le m$, the computed corrections, we can define 
approximation spaces $\Wc_{m}^k = span\{w_i^{(k)}\}_{i=1}^m\subset \Sc_{n_k}^k$ (with dimension at most $m$), and look for an approximation of the form $u_m = \sum_{i_1=1}^m \hdots \sum_{i_r=1}^m \alpha_{i_1\hdots i_r} \otimes_{k=1}^r w_{i_k}^{(k)} \in \otimes_{k=1}^r \Wc^k_m$ (namely Tucker tensor format). 
The update problem then consists in solving
  \begin{align}
\min_{\alphab=(\alpha_{i_1\hdots i_r}) \in \Rbb^{m\times\ldots\times m}}   \Vert u-\sum_{i_1,\hdots,i_r} \alpha_{i_1\hdots i_r}  \otimes_{k=1}^r w_{i_k}^{(k)}\Vert_Q^2 + 
\lambda' \Vert \alphab \Vert_1, \label{eq:algo_update_regularized_tucker}
\end{align}
where $\Vert \alphab\Vert_1=\sum_{i_1 ,\hdots, i_r}\lvert \alpha_{i_1 \hdots  i_r}\lvert$.
This updating strategy can yield significant improvements of convergence. However,
it is clearly unpractical for high dimension $r$ since the dimension $m^r$ of the representation grows exponentially with $r$. For high dimension, other types of representations should be introduced, such as hierarchical tensor representations.
\end{remark}

Algorithm \ref{alg:sparse_rank_m} details the updated greedy construction of sparse low rank approximations.

\begin{algorithm}
\caption{Updated greedy algorithm for sparse low rank approximation of a function $u$.}
\label{alg:sparse_rank_m}
\begin{algorithmic}[1]
\REQUIRE vector of evaluations $\z = (u(y^1),\hdots,u(y^Q))^T \in\Rbb^Q$ and maximal rank $M$.
\ENSURE Sequence of approximations $u_m = \sum_{i=1}^m \alpha_i w_i$, where  $w_i \in \Rc_1$ and $\alphab=(\alpha_1,\hdots,\alpha_m)\in\Rbb^m$.
\STATE Set $u_0=0$
\FOR {$m=1,\ldots,M$}
\STATE Evaluate the vector $\z_{m-1}=(u_{m-1}(y^1),\hdots,u_{m-1}(y^Q))^T\in\Rbb^Q$
\STATE \label{algo3-correction} Compute a sparse rank-one approximation $w_m=\otimes_{k=1}^r w_m^{(k)}$ using Algorithm \ref{sparse_rank_one_algo} for input vector of evaluations $\z-\z_{m-1}$. 
\STATE Evaluate matrix $\Wb\in\Rbb^{Q\times m}$ with components $(\Wb)_{qi}=w_i(y^q)$.
\STATE   \label{algo3-update} Compute $\alphab\in \Rbb^m$ with Algorithm \ref{alg:loo_cv} for input vector $\z\in\Rbb^Q$ and matrix $\Wb\in\Rbb^{Q\times m}$.
\ENDFOR
\end{algorithmic}
\end{algorithm}
 
 \subsubsection{Rank Selection}\label{sec:rankselection}
Let us note that Algorithm \ref{alg:sparse_rank_m}
 not only generates an approximation of the function but a sequence of approximations $\{u_m\}_{m =1}^M$ with increasing canonical ranks. To select the best rank, we use a $k$-fold cross validation method. The overall procedure is as follows:  
\begin{itemize}
\item Split sample set $S=\{1,\hdots,Q\}$ into $k$ disjoint subsamples $\{V_i\}_{i=1}^k$, $V_i\subset S$, of approximately the same size, and let  $S_i=S\setminus V_i$.
\item For each subsample, 
run Algorithm \ref{alg:sparse_rank_m} from the reduced vector of evaluations $\z_{S_i}$ (test set) 
to obtain the sequence of model approximations $\{u_m^{S_i}\}_{1\leq m\leq M}$.
Compute the corresponding mean squared errors $\{\varepsilon_1^{S_i},\ldots,\varepsilon_M^{S_i}\}$ 
from the validation set of evaluations $\z_{V_i}$.
\item For $m=1,\ldots,M$, compute the $k$-fold cross validation error $\bar{\varepsilon}_m=\frac{1}{k} \sum_{i=1}^k \varepsilon_m^{S_i}$.
\item Select optimal rank $m_{op}=\text{argmin}_{1\leq m\leq M} \bar{\varepsilon}_m$.
\item Run Algorithm \ref{alg:sparse_rank_m} with the whole data set $\z$ for computing $u_{m_{op}}$. 
\end{itemize}

\bigskip
\begin{remark}[Non greedy approximation]\label{r1_closed}
{\color{black}
A non greedy variant of Algorithm \ref{alg:sparse_rank_m} for computing a rank-$m$ approximation would consist in performing approximation directly
in $\Rc_m$ (as proposed in \cite{DOO12} without sparsity-inducing regularization). This strategy usually yields better approximations for the same rank $m$. However, the set $\Rc_m$ for $m\geq 2$ being not closed (see \cite{HAC12} lemma 9.11 pg 255), minimization  in $\Rc_m$ is an ill-posed problem which usually requires regularization. The proposed construction based on successive rank-one corrections yields a suboptimal low-rank decomposition but it has several advantages: successive minimization problems in $\Rc_1$ are well-posed, it requires the learning of a small number of parameters at each iteration. 
}
\end{remark}

\section{Application examples}\label{applex}

{\color{black}In this section, we validate the proposed algorithm on several benchmark problems. The purpose of the first example on Friedman function in section \ref{sec:friedman} is to highlight the benefit singly of the greedy low rank approximation by giving some hints on the number of samples needed for a stable approximation. 
The three following examples then exhibit the robustness of the $l_1$-regularization within the low rank approximation: 
\begin{itemize}
	\item by correctly detecting sparsity when appropriate approximation space is introduced as in the checker board function case presented in section \ref{sec:checkerboard},
	\item or just by looking for the simplest representation with respect to the number of samples as in the examples of sections \ref{sec:rastrigin} and \ref{sec:twoplates}.
\end{itemize}
}

In all the examples, for the purpose of estimating the approximation errors, we introduce the relative error $\varepsilon(u_m,u)$ between the function $u$ and a rank-$m$ approximation $u_m$, estimated with Monte Carlo integration with $Q'=1000$ samples:
\begin{align}
\varepsilon (u_m,u)=\frac{\norm{ u_m-u}_{Q'}}{\norm{ u }_{Q'}}.
\end{align}

{\color{black}Let us also define the sparsity ratio $\varrho_1$ of a rank-$1$ approximation $w = \otimes_{k=1}^d w^{(k)}$ as:
$$
\varrho_1 =\frac{\sum_{k=1}^d\beta(w^{(k)})}{\sum_{k=1}^d {n_k}}
$$
where $\beta(w^{(k)})=\norm{\w^{(k)}}_0$ gives the number of non zero coefficients in the vector of coefficients $\w^{(k)}$. In short, $\varrho_1$ is the ratio of total number of non zero parameters to the total number of parameters in rank one tensor representation. We also define the total sparsity ratio $\varrho_{m}$ and the partial sparsity ratio $\varrho_{m}^{(k)}$ in the dimension $k$ of a rank-$m$ approximation $u_m = \sum_{i=1}^m \otimes_{k=1}^d w_i^{(k)}$ as: 
$$
\varrho_{m} =\frac{\sum_{i=1}^m\sum_{k=1}^d\beta(w^{(k)}_i)}{m\times \sum_{k=1}^d{n_k}} \quad \text{and } \quad \varrho_{m}^{(k)} =\frac{\sum_{i=1}^m\beta(w^{(k)}_i)}{m\times {n_k}}.
$$

}

{\color{black}
\subsection{Analytical model: Friedman function}\label{sec:friedman}
Let us consider a simple benchmark problem namely the Friedman function of dimension $d=5$ also considered in \cite{BEL11}:
$$
u(\xib)=10\sin(\pi \xi_1 \xi_2)+20(\xi_3-0.5)^2+10\xi_4+5\xi_5,
$$
 where $\xi_i, i=1,\ldots,5$, are uniform random variables over [0,1].
In this section, we consider low rank tensor subsets without any sparsity constraint, that means we do not perform $\ell_1$ regularization in the alternating minimization algorithm. The aim is to estimate numerically the number of function evaluations necessary in order to obtain a stable approximation in low rank tensor format. 

The analysis on the number of function evaluations is here based upon results of \cite{NOB11} where it is proven that, for a stable approximation of a monovariate function with optimal convergence rate using ordinary least squares on polynomial spaces, the number of random sample evaluations required scales quadratically with the dimension of the polynomial space, i.e. $Q\sim(p+1)^2$ where $p$ is the maximal polynomial degree. 
This result is supported with numerical tests on monovariate functions and on multivariate functions  
that show that indeed choosing $Q$ scaling quadratically with the dimension of the polynomial space $P$ is robust 
 while $Q$ scaling linearly with $P$ may lead to an ill conditioned problem and hence an unstable approximation, depending on $P$ and the dimension $d$. 

The following numerical tests aim at bringing out a similar type of rule for choosing the number of samples $Q$ for a low rank approximation of a multivariate function constructed in a greedy fashion according to Algorithm \ref{alg:sparse_rank_m}, but with ordinary least squares in step \ref{step:alg2_l1} of Algorithm \ref{algo2-loop}, and given an isotropic tensor product polynomial approximation space with maximum degree $p$ in all dimensions. We first consider a rank one approximation of the function in $\Pbb_{p}\otimes \cdots \otimes \Pbb_{p}$ where $\Pbb_{p}$ denotes the polynomial approximation space of maximal degree $p$ in each dimension from $1$ to $d$.
Given the features above and considering the algorithm for the construction of the rank one element of order $d$, we consider the following  rule:
$$Q=cd(p+1)^\alpha$$
where $c$ is a positive constant and $\alpha=1$ (linear rule) or $2$ (quadratic rule).  
In the following analyses of the current section, we plot the mean $\varepsilon(u_m,u)$ over 51 sample set repetitions in order to eliminate any dependence on the sample set of a given size.
In figure \ref{friedman1_rank_one}, we compare the error of rank one approximation $\varepsilon(u_1,u)$ with respect to the Legendre polynomial degree $p$ using both linear rule (left) and quadratic rule (right) for different values of $c$ (ranging from 1 to 20 in the linear rule and 0.5 to 3 in quadratic rule). As could have been expected, we find that the linear rule yields a deterioration for small values of $c$ whereas the quadratic rule gives a stable approximation with polynomial degree.

\begin{figure}[h!]\centering
\subfigure[]{\includegraphics[width=0.45\linewidth]{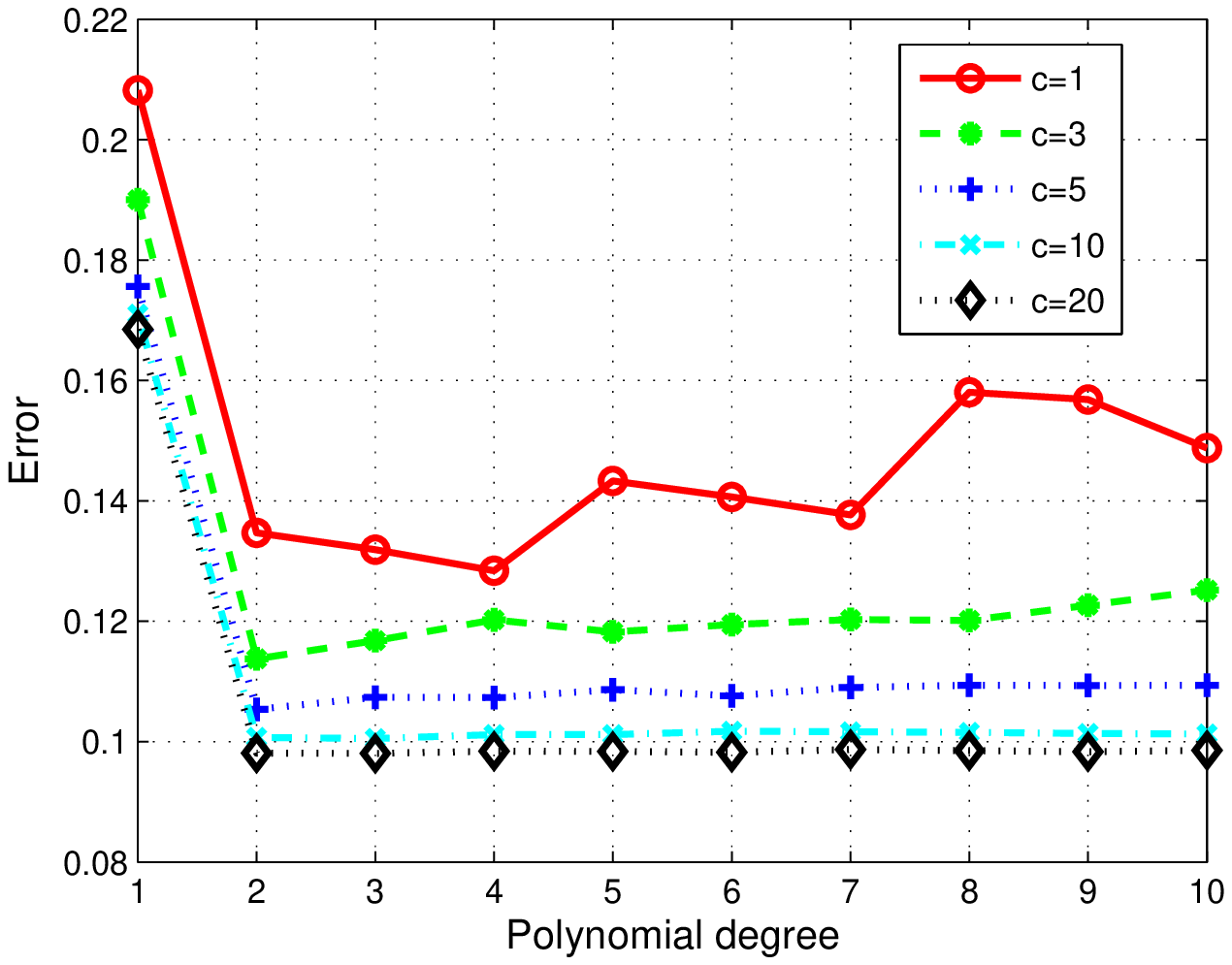}}
\subfigure[]{\includegraphics[width=0.45\linewidth]{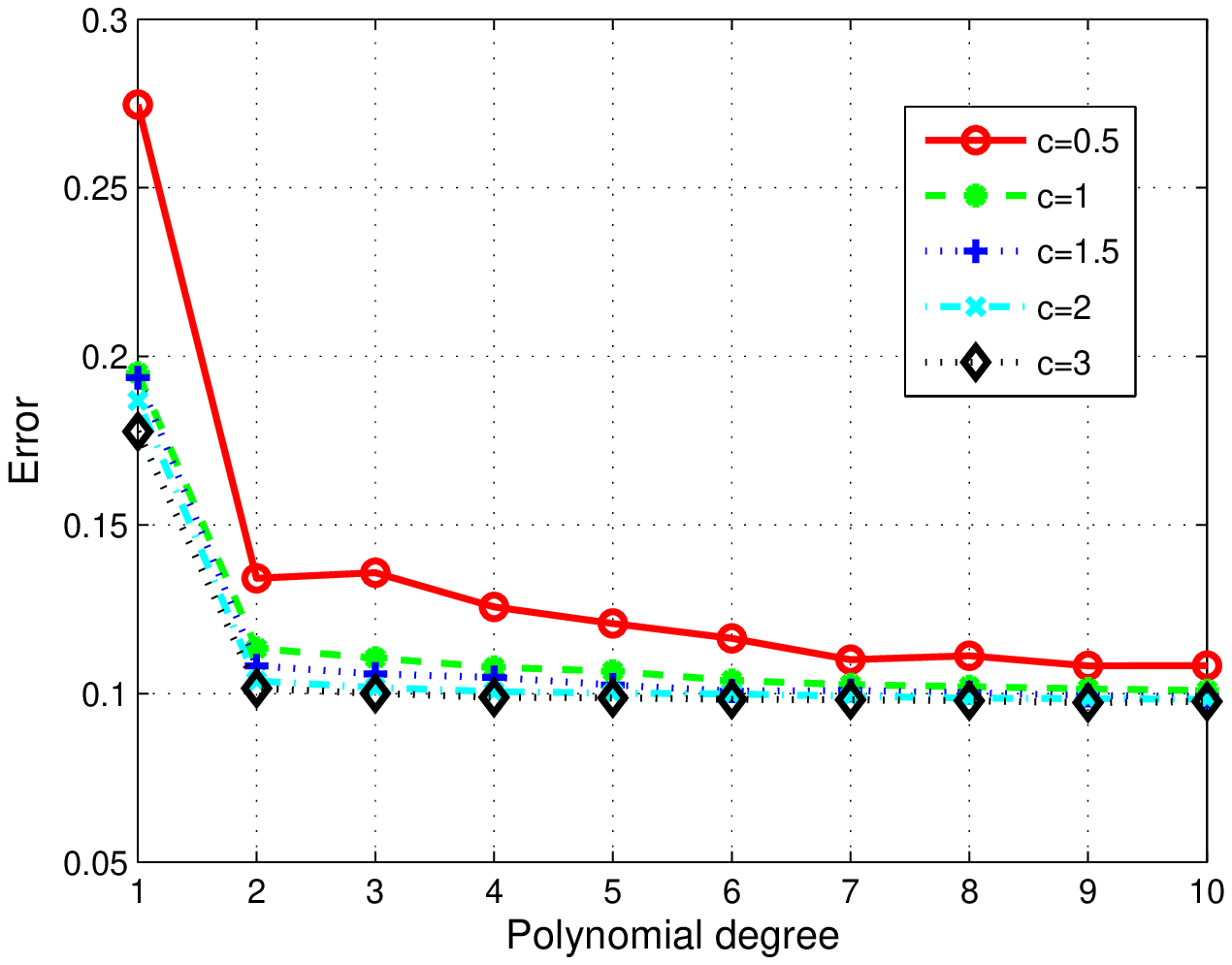}}
  \caption{Evolution of rank one approximation error $\varepsilon(u_1,u)$ with respect to polynomial degree $p$ with (a) $ Q=cd(p+1)$ and (b) $Q=cd(p+1)^2$ samples with several values of $c$.}
  \label{friedman1_rank_one}
\end{figure}

For higher rank approximations, the total number of samples needed will have a dependence on the approximation rank $m$. Thus we modify sample size estimates and consider the rule $$Q= cdm(p+1)^\alpha$$ 
with $\alpha=1$ (linear rule) or $2$ (quadratic rule). 
In figure \ref{friedman1_high_rank}, we plot approximation error  $\varepsilon(u_m,u)$ using linear rule (left) and quadratic rule (right) for $m=2,3,4$ and different values of $c$. We find again that quadratic rule gives a stable approximation for $c\geq1$.

\begin{figure}[h!]\centering
\subfigure[$m=2$]{\includegraphics[width=0.45\linewidth]{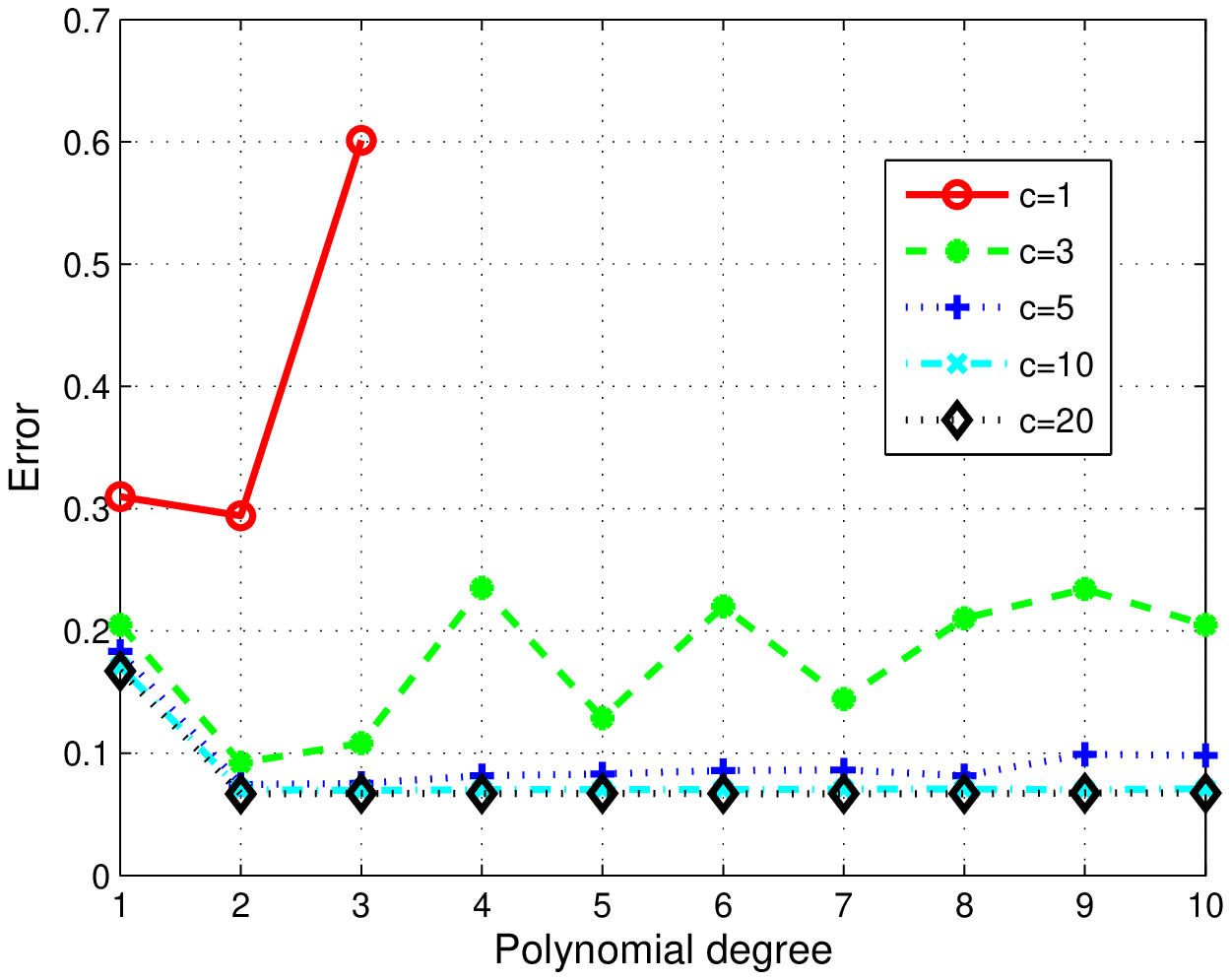}}
\subfigure[$m=2$]{\includegraphics[width=0.45\linewidth]{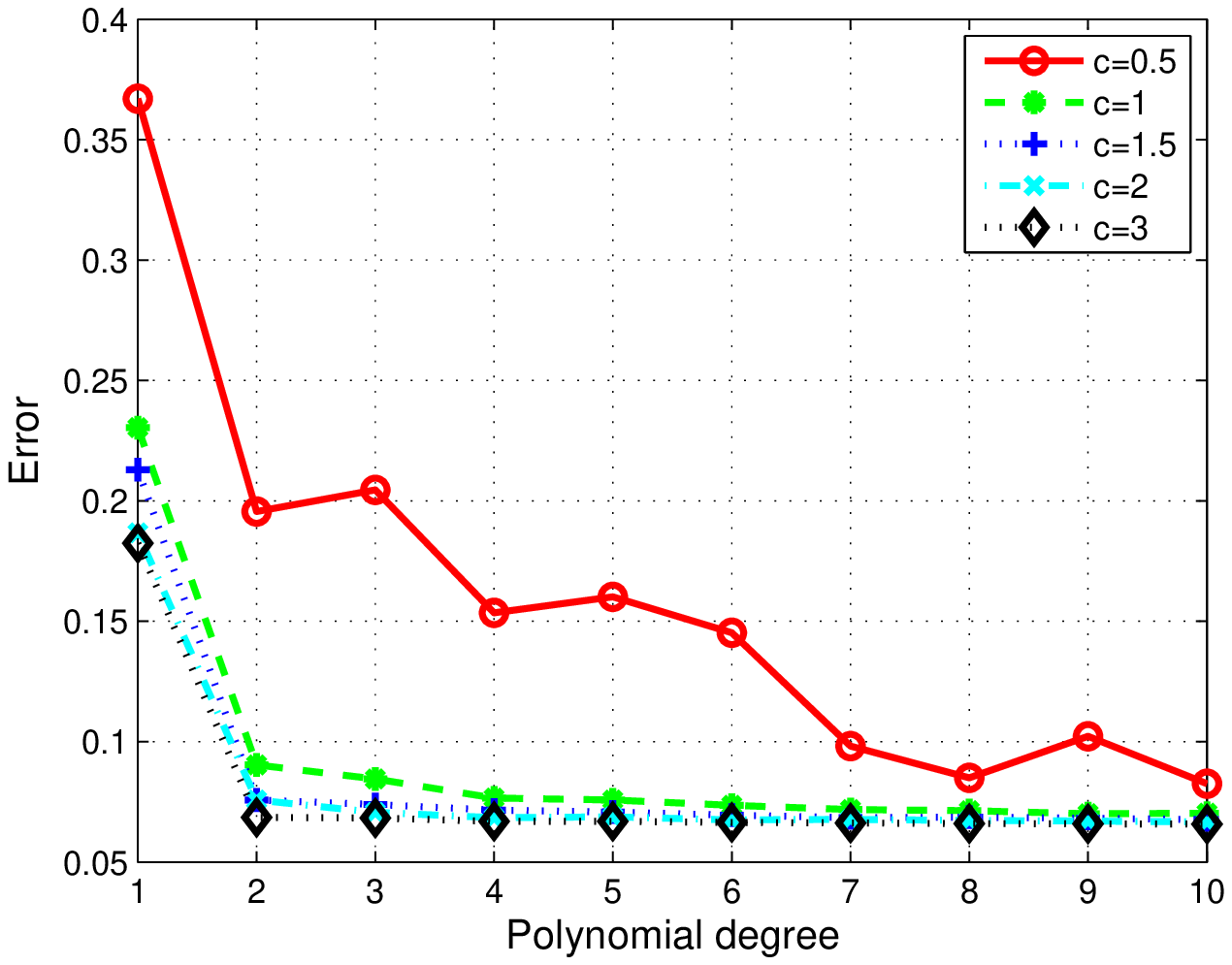}}
\\
\subfigure[$m=3$]{\includegraphics[width=0.45\linewidth]{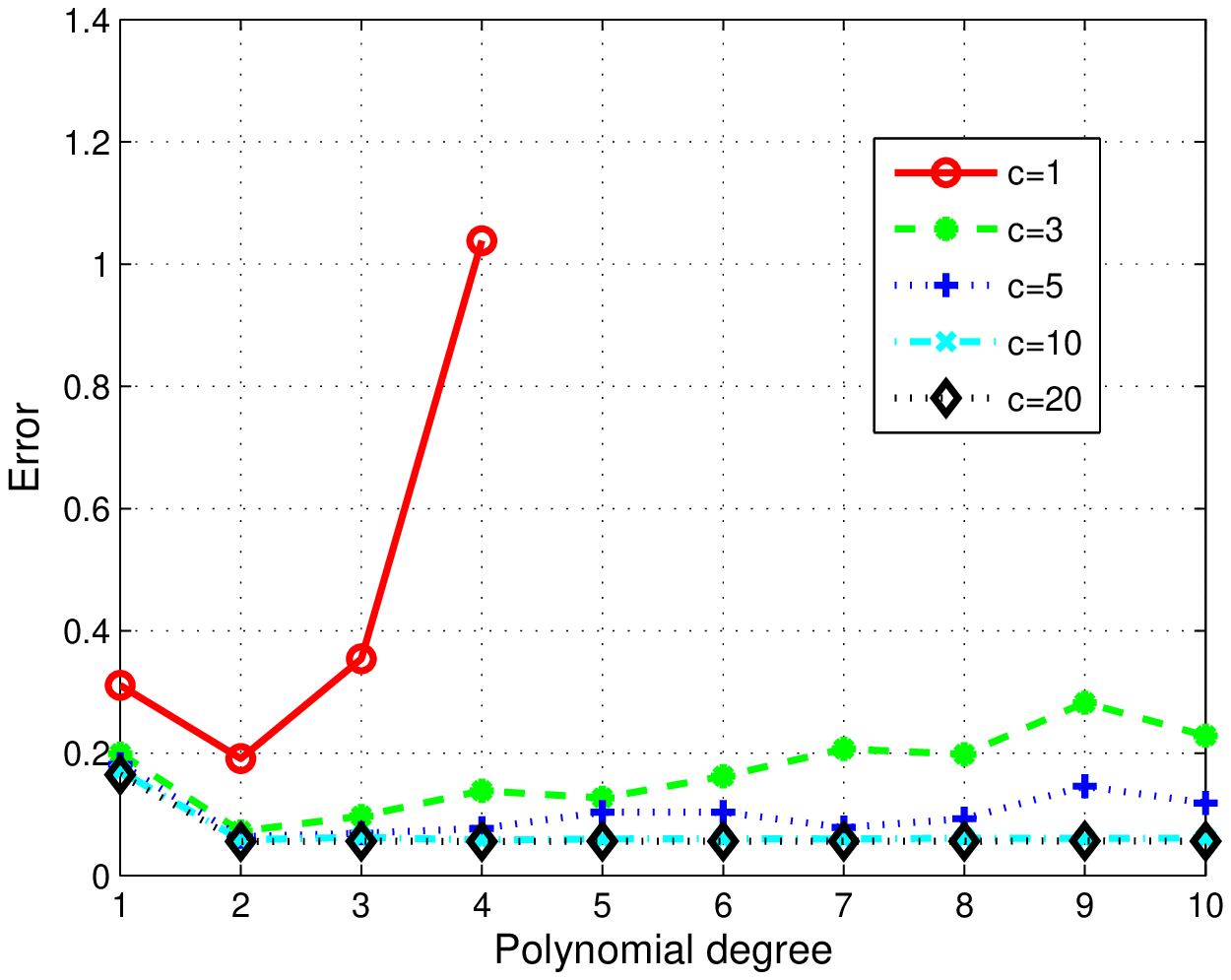}}
\subfigure[$m=3$]{\includegraphics[width=0.45\linewidth]{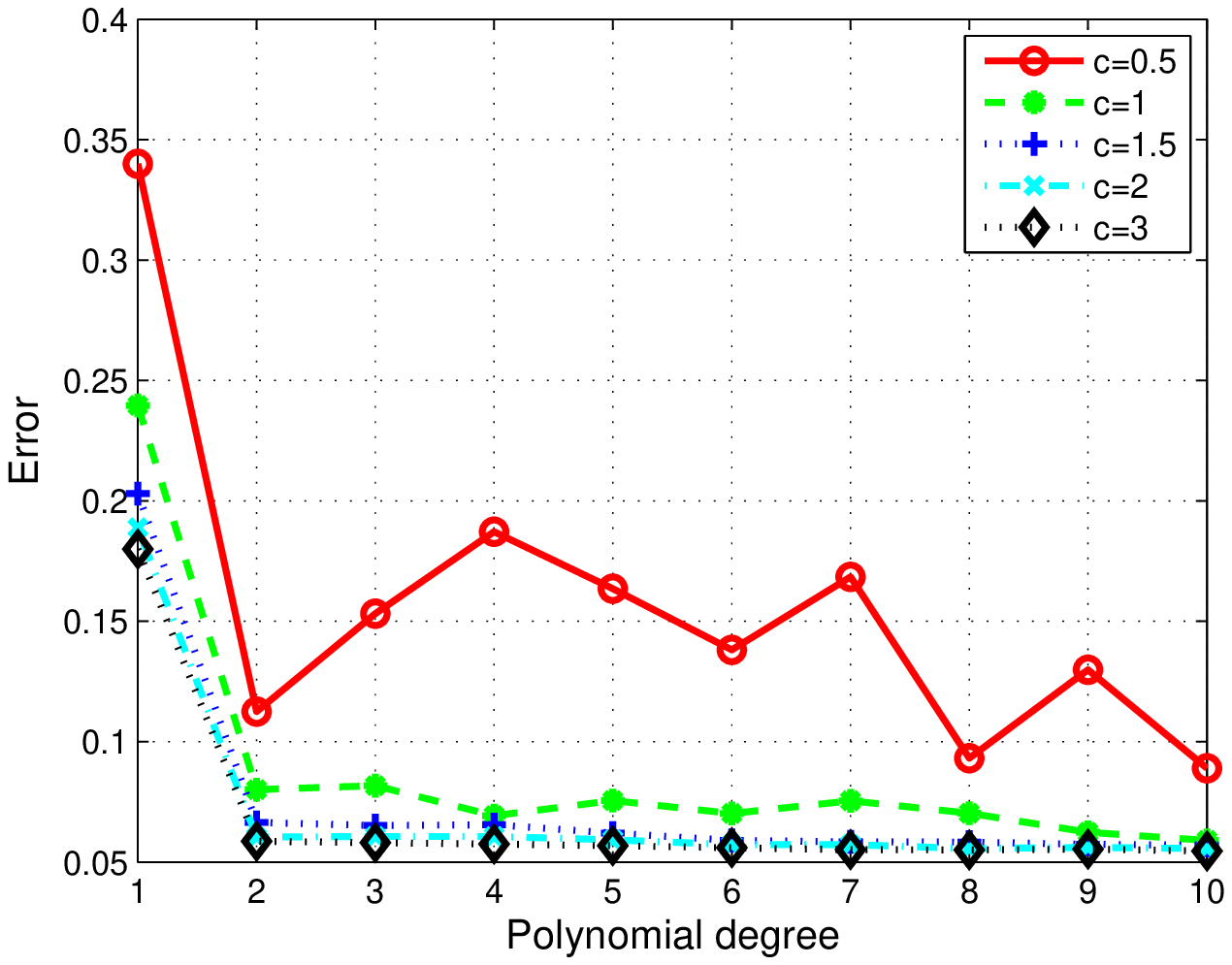}}
\\
\subfigure[$m=4$]{\includegraphics[width=0.45\linewidth]{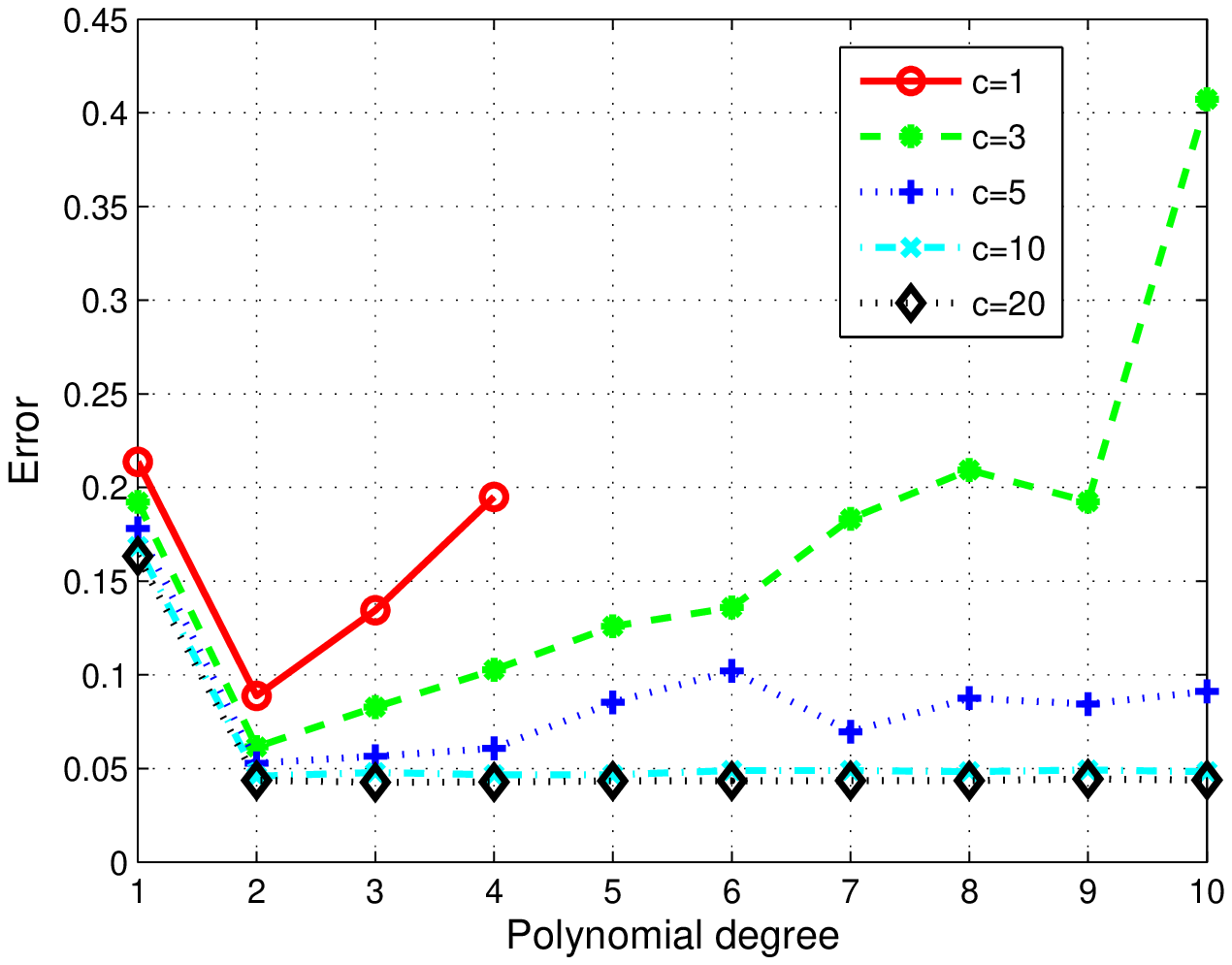}}
\subfigure[$m=4$]{\includegraphics[width=0.45\linewidth]{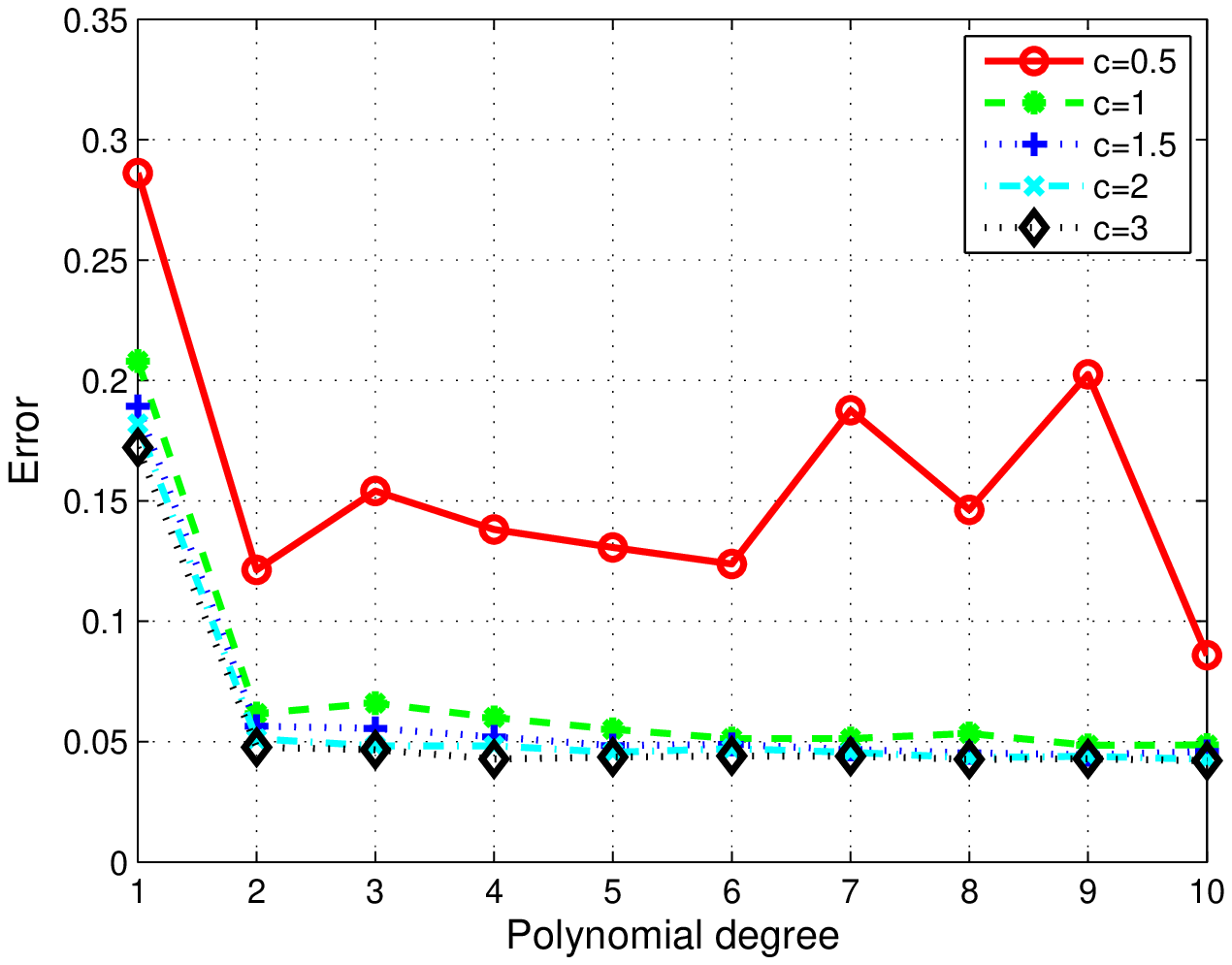}}
  \caption{Evolution of approximation error $\varepsilon(u_m,u)$ with respect to polynomial degree $p$ with $ Q=cd(p+1)$ (left column) and $Q=cd(p+1)^2$ (right column)  with several values of $c$ and for $m=2,3,4$.}
  \label{friedman1_high_rank}
\end{figure}

In order to analyze the accuracy of the rank-$m$ approximation with respect to $m$, in figure \ref{friedman1_rank} we plot the error $\varepsilon(u_m,u)$ with respect to the polynomial degree $p$ for different values of $m$ using $Q=dm(p+1)^2$, that is $c=1$. As the number of samples $Q$ increases with rank $m$ using this rule, more information on the function is given enabling for higher rank approximations to better represent the possible local features of the solution. We thus find that the approximation is more accurate as the rank increases.
\begin{figure}[h!]
\centering
\includegraphics*[width=0.5\textwidth]{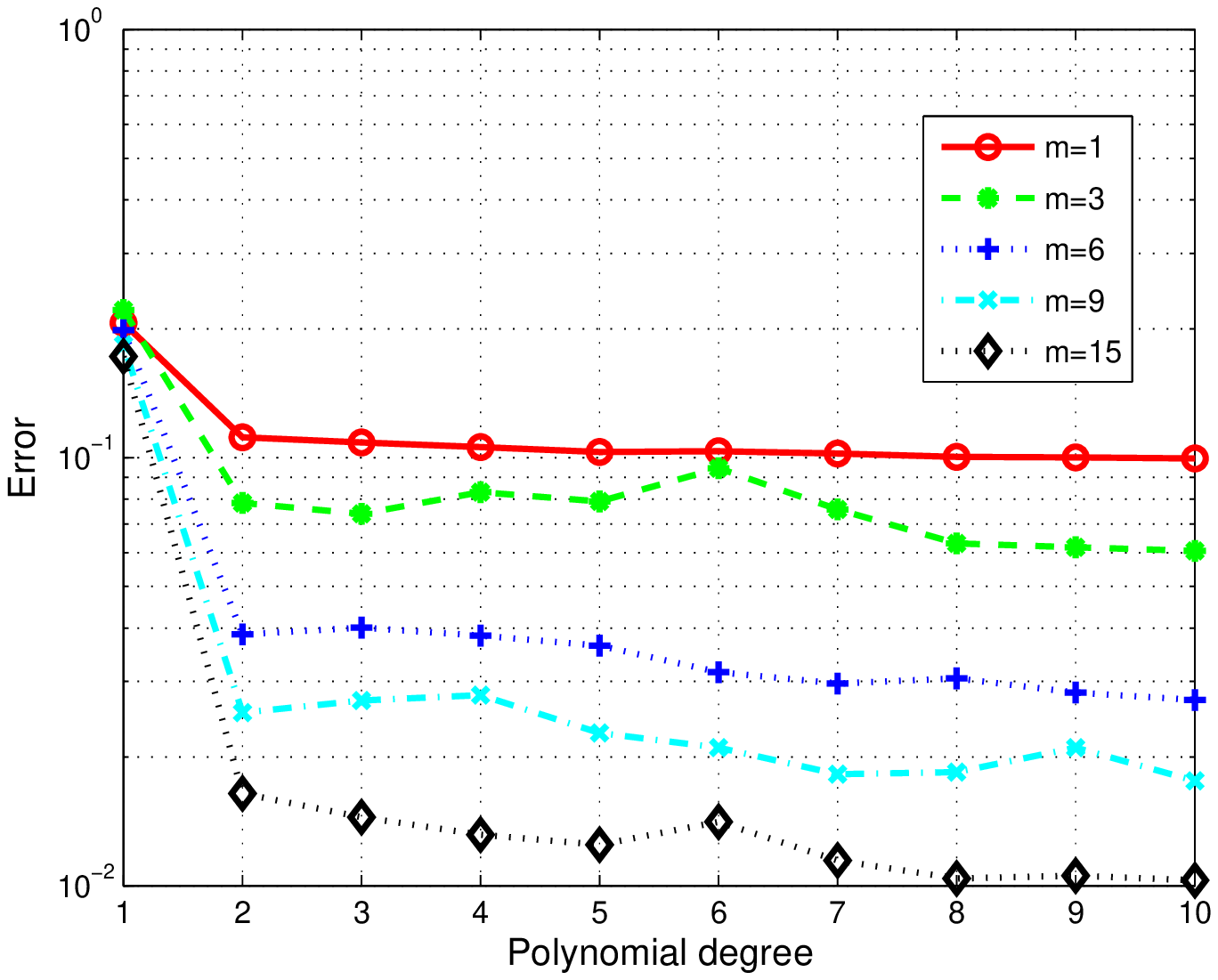}
\caption{Evolution of error $\varepsilon(u_m,u)$ with respect to polynomial degree $p$ for different values of $m$ using sample size given by quadratic rule $Q= dm(p+1)^2$.}
\label{friedman1_rank}
\end{figure}

From this example, we can draw the following conclusions:
\begin{itemize}
\item a heuristic rule to determine the number of samples needed in order to have a stable low rank approximation grows only linearly with dimension $d$ and rank $m$ and is given by $Q=dm(p+1)^2$,
\item better solutions are obtained with high rank approximations, provided that enough model evaluations are available.
\end{itemize} 

Quite often in practice, we do not have enough model evaluations and hence we may not be able to achieve good approximations with limited sample size. This is particularly true for certain classes of non smooth functions. One possible solution is a good choice of bases that are sufficiently rich (such as piecewise polynomials or wavelets) and that can capture simultaneously both global and local features of the model function. However, the sample size may not be sufficient enough to obtain good approximations with ordinary least squares in progressive rank one corrections due to large number of basis functions. We illustrate in section \ref{sec:checkerboard} and \ref{sec:rastrigin} that, in such cases, performing approximation in $sparse$ low rank tensor subsets (i.e. using $\ell_1$ regularization in alternating minimization algorithm) allows more accurate approximation of the model function. In addition, we illustrate in section \ref{sec:twoplates} that approximation in sparse low rank tensor subsets leads to a relatively stable approximation with limited number of samples even for high degree polynomial spaces.

}

\subsection{Analytical model: Checker-board function}\label{sec:checkerboard}

\subsubsection{Function and approximation spaces}
\begin{figure}[h!]\centering
\includegraphics[width=0.4\linewidth]{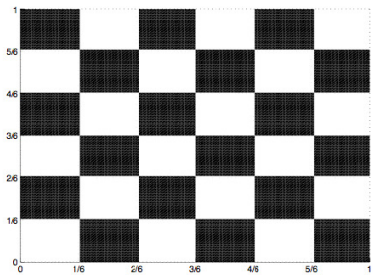}
\caption{Checker-board function.
}
 \label{fig:checker_board}
\end{figure}

We now test Algorithm \ref{alg:sparse_rank_m} on the so-called checker-board function $u(\xi_1,\xi_2)$ of dimension $d=2$ illustrated in figure \ref{fig:checker_board}. {\color{black}The purpose of this test is to illustrate that, given appropriate bases, in this case piecewise polynomials, Algorithm \ref{alg:sparse_rank_m} allows the detection of sparsity and hence construction of a sequence of optimal sparse rank-$m$ approximations with few samples.}

Random variables $\xi_1$ and $\xi_2$
are independent and uniformly distributed on $[0,1]$. The checker-board function is a rank-$2$ function 
$$u(\xi_1,\xi_2)=\sum_{i=1}^2 w_i^{(1)}(\xi_1)w_i^{(2)}(\xi_2)$$
with $w_1^{(1)}(\xi_1)=c(\xi_1)$, $w_1^{(2)}(\xi_2)=1-c(\xi_2)$, $w_2^{(1)}(\xi_1)=1-c(\xi_1)$ and $w_2^{(2)}(\xi_2)=c(\xi_2)$  where $c(\xi_k)$ is the crenel function defined by:
\begin{align*}
c(\xi_k)=\begin{cases}
1 \quad \text{ on } [0,\frac{1}{6}[+2n\frac{1}{6}, \,n=0,1,2 \\
0 \quad \text{ on } [\frac{1}{6},\frac{2}{6}[+2n\frac{1}{6}, \,n=0,1,2 \\
\end{cases}.
\end{align*}


For approximation spaces $\Sc^k_{n_k}$, $k\in\{1,2\}$, we introduce piecewise polynomials of degree $p$ defined on a uniform partition of $\Xi_k$ composed by $s$ intervals, corresponding to $n_k=s(p+1)$. We denote by $\Sc^k_{n_k} = \Pbb_{p,s}$ the corresponding space (for ex. $\Pbb_{2,3}$ denotes piecewise polynomials of degree 2 defined on the partition $\{(0,\frac{1}{3}),(\frac{1}{3},\frac{2}{3}),(\frac{2}{3},1)\}$). We use an orthonormal basis composed of functions whose supports are one element of the partition and whose restrictions on these supports are rescaled Legendre polynomials. 

Note that when using a partition into $s=6n$ intervals, $n\in\Nbb$, then the checker-board function can be exactly represented, that means $u\in \Pbb_{p,6n}\otimes \Pbb_{p,6n}$ for all $p$ and $n$. Also, the solution admits a sparse representation in $\Pbb_{p,6n}\otimes \Pbb_{p,6n}$ since an exact representation is obtained by only using piecewise constant basis functions ($u\in \Pbb_{0,6n}\otimes \Pbb_{0,6n}$). The effective dimensionality of the checker-board function is $2\times 2\times 6 = 24$, which corresponds to the number of coefficients required for storing the rank-2 representation of the function. 
We expect our algorithm to detect the low-rank of the function and also to detect its sparsity. 

\subsubsection{Results}
Algorithm \ref{alg:sparse_rank_m} allows the construction of a sequence of sparse rank-$m$ approximations $u_m$ in $\Sc^1_{n_1}\otimes \Sc^2_{n_2}$. We estimate optimal rank-{$m_{op}$} using 3-fold cross validation (see Section \ref{sec:rankselection}).

In order to illustrate the accuracy  of approximations in sparse low rank tensor subsets, we compare the performance of $\ell_1$-regularization within the alternating minimization algorithm (step \ref{algo3-correction} of Algorithm \ref{alg:sparse_rank_m}) with no regularization (OLS) and the $\ell_2$-regularization (see Remark \ref{algo_versions} for the description of these alternatives). 
Table \ref{tab:1} shows the error $\varepsilon(u_{m_{op}},u)$ obtained for the selected optimal rank $m_{op}$, without and with updating step \ref{algo3-update} of Algorithm \ref{alg:sparse_rank_m}, and for the different types of regularization during the correction step \ref{algo3-correction} of Algorithm \ref{alg:sparse_rank_m}. 
The results are presented for a sample size $Q=200$ and for different function spaces $\Sc^1_{n_1}= \Sc^2_{n_2}=\Pbb_{p,s}$.
$P$ denotes the dimension of the space $\Sc^1_{n_1}\otimes \Sc^2_{n_2}$. 
We observe that, for $\Pbb_{2,3}$, the solution is not sparse on the corresponding basis and 
$\ell_1$-regularization does not provide a  better solution than $\ell_2$-regularization since the approximation space is not adapted. However, when we choose 
function spaces that are sufficiently rich for the solution to be sparse, we see that $\ell_1$-regularization within the alternating minimization algorithm
outperforms other types of regularization and yields low rank approximations of the function almost at the machine precision. {\color{black}This is because $\ell_1$-regularization is able to select non zero coefficients corresponding to appropriate basis functions of the piecewise polynomial approximation space. For instance, when $\Pbb_{5,6}$ is used as the approximation space, only 3 (out of 36) non zero coefficients corresponding to piecewise constant bases are selected by $\ell_1$ regularization along each dimension in each rank one element (that is the sparsity ratio is $\varrho_1\approx 0.08$ for each rank-one element), thus yielding an almost exact recovery}. We also find that $\ell_1$-regularization allows recovering the exact rank-$2$ approximation of the function.

\begin{table}[h!]
\footnotesize
\caption{Relative error $\varepsilon(u_{m_{op}},u)$ and optimal rank $m_{op}$ estimation of Checker-board function with various regularizations for $Q=200$ samples. $P$ is the dimension of the approximation space. ({`-' indicates that none of the rank-one elements were selected during the update step}).}
\label{tab:1}
\begin{center}
\begin{tabular}{l|llll|cccc|cccc}
\toprule
& \multicolumn{4}{c|}{Ordinary Least Square} & \multicolumn{4}{c|}{$\ell_2$} & \multicolumn{4}{c}{$\ell_1$}\\
\hline 
&\multicolumn{2}{c}{No update} & \multicolumn{2}{c|}{Update} & \multicolumn{2}{c}{No update} & \multicolumn{2}{c|}{Update} &\multicolumn{2}{c}{No update} & \multicolumn{2}{c}{Update}\\
\hline 
Approximation space & Error & $m_{op}$ & Error & $m_{op}$ & Error & $m_{op}$ & Error & $m_{op}$ & Error & $m_{op}$ & Error & $m_{op}$\\
\hline 
$\Rc_m(\Pbb_{2,3}\otimes\Pbb_{2,3}),P=9^2$ & 0.527 & 2 & 0.527 & 2 & 0.508 & 2 & 0.508 & 2 & 0.507 & 2 & 0.507 & 2\\
$\Rc_m(\Pbb_{2,6}\otimes \Pbb_{2,6}),P=18^2$ & 0.664 & 2 & 0.664 & 2 & 0.061 & 8 & 0.061 & 8 &$1.96\, 10^{-12}$ & 4 & $2.41\, 10^{-13}$ & 2\\

$\Rc_m(\Pbb_{2,12}\otimes \Pbb_{2,12}),P=36^2$& 20.92 & 1 & - & - & 0.568 & 10 & 0.566 & 4 & $1.93\, 10^{-12}$ & 2 & $1.1\, 10^{-12}$ & 2\\

$\Rc_m(\Pbb_{5,6}\otimes \Pbb_{5,6}),P=36^2$& 31.27 & 1 & - & - & 0.624 & 10 & 0.623 & 3 & $1.22\, 10^{-12}$ & 2 & $7.93\, 10^{-13}$ & 2\\ 

$\Rc_m(\Pbb_{10,6}\otimes \Pbb_{10,6}),P=66^2$& 9648.8 & 1 & - & - & 0.855 & 10 & 0.855 & 10 & $1.21\,10^{-12}$ & 2 & $7.88\, 10^{-13}$ &2 \\

\bottomrule
\end{tabular}
\end{center}
\end{table}
From this analytical example, several conclusions can be drawn:
\begin{itemize}
\item $\ell_1$-regularization in alternating least squares algorithm is able to detect sparsity and hence gives very accurate approximations using few samples as compared to OLS and $\ell_2$-regularizations,
\item updating step selects the most pertinent rank-one elements and gives an approximation of the function with a lower effective rank.
\end{itemize} 

\subsection{Analytical model: Rastrigin function}\label{sec:rastrigin}


{\color{black} For certain classes of non smooth functions, wavelet bases form an appropriate choice as they allow the simultaneous description of global and local features \cite{LEM04}.} In this example, we demonstrate the use of our algorithm with polynomial wavelet bases by studying 2-dimensional Rastrigin function. The function is given by

$$
u(\xib)=20+\sum_{i=1}^2(\xi_i^2-10\cos(2\pi\xi_i))
$$ 
where $\xi_1,\xi_2$ are independent random  variables uniformly distributed in [-4,4]. 
%
%

We consider two types of approximation spaces $\Sc^k_{n_k}$, $k\in\{1,2\}$:
\begin{itemize}
\item spaces of polynomials of degree 7, using Legendre polynomial chaos basis, denoted $\Pbb_{7} $,
\item spaces of polynomial wavelets with degree 4 and resolution level 3, denoted $\mathbb W_{4,3}$.
\end{itemize}

\begin{figure}[h!]
\centering
\subfigure[]{\includegraphics[width=7.5cm]{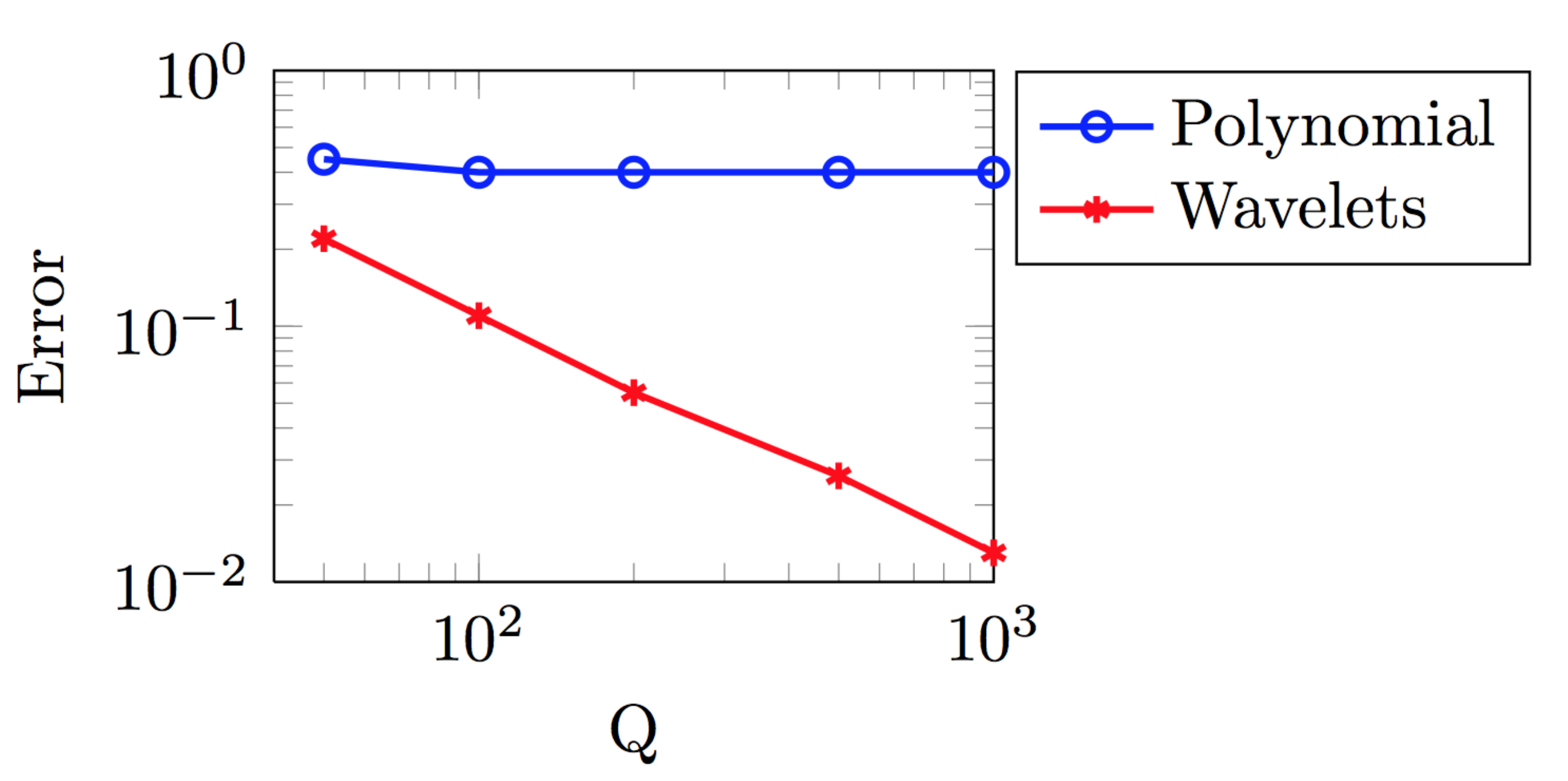}}
\subfigure[]{\includegraphics[width=7.cm]{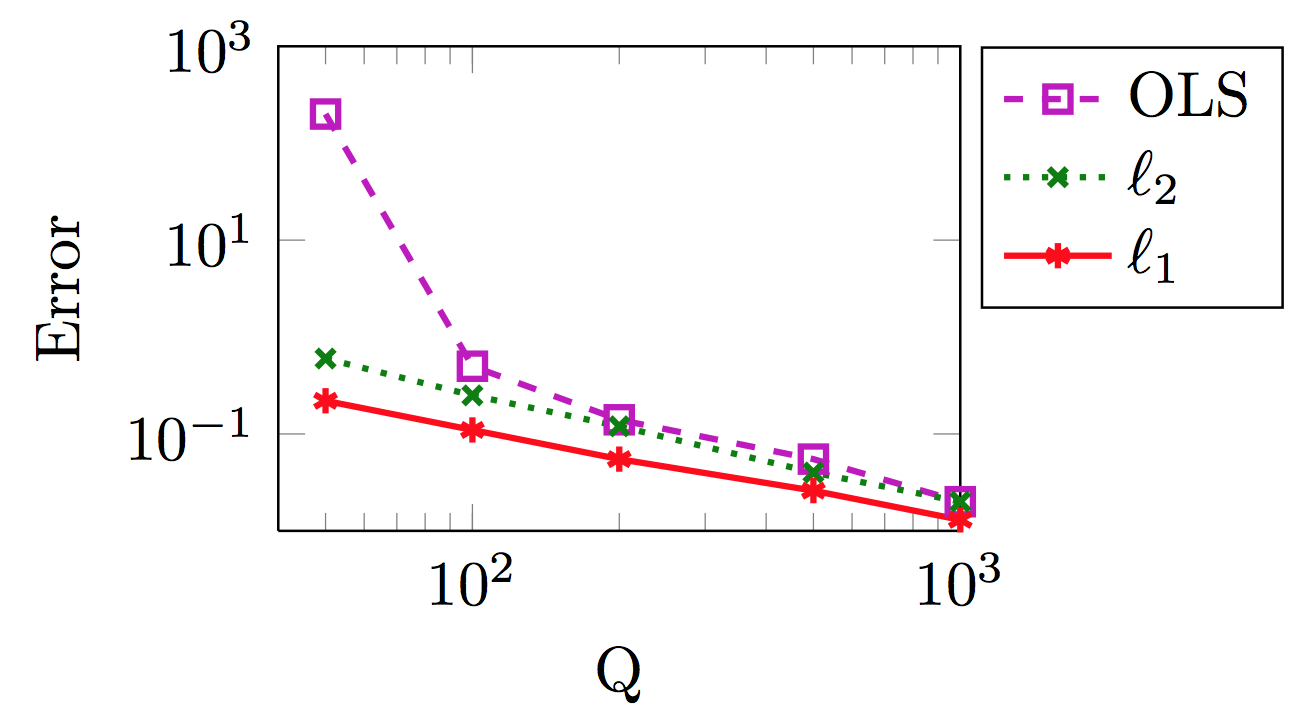}}
\caption{{\color{black}Evolution of error $\varepsilon(u_{m_{op}},u)$ with respect to the number of samples $Q$. Approximations obtained with Algorithm \ref{alg:sparse_rank_m} with optimal rank selection: (a) for the two different approximation spaces  $\Pbb_{7}\otimes \Pbb_{7}$ ($P=64$) and $\mathbb W_{4,3}\otimes \mathbb W_{4,3}$ ($P=1600$) and (b) for different types of regularizations with approximation space $\mathbb W_{4,3}\otimes \mathbb W_{4,3}$}.}
\label{rastrigin_review_1}
\end{figure}

We compute a sequence of sparse rank-$m$ approximations $u_m$ in $\Sc_{n_1}^1 \otimes \Sc_{n_2}^2$ using Algorithm \ref{alg:sparse_rank_m} and an optimal rank approximation $u_{m_{op}}$ is selected using 3-fold cross validation (see the rank selection strategy in section \ref{sec:rankselection}). 
Figure \ref{rastrigin_review_1}(a) shows the convergence of this optimal approximation with respect to the sample size $Q$ for the two different approximation spaces. 
We find that the solution obtained with classical polynomial basis functions is inaccurate and does not improve with increase in sample size. Thus, polynomial basis functions are not a good choice to obtain a reasonably accurate estimate. On the other hand, when we use wavelet approximation bases, the approximation error reduces progressively with increase in sample size. {\color{black} Figure \ref{rastrigin_review_1}(b) shows the convergence of the optimal wavelet approximation with respect to the sample size $Q$ for different regularizations within the alternated minimization algorithm of the correction step. The $\ell_1$ regularization is more accurate when compared to both OLS and $\ell_2$ regularization, particularly for few model evaluations. We can thus conclude that a good choice of basis functions is important in order to fully realize the potential of sparse $\ell_1$ regularization in the tensor approximation algorithm.}

\begin{figure}[h!]
\centering
\includegraphics*[width=0.50\textwidth]{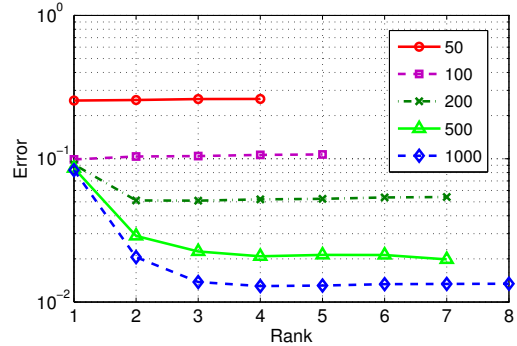}
\caption{Evolution of error $\varepsilon(u_{m_{op}},u)$ of approximations obtained using Algorithm \ref{alg:sparse_rank_m} with respect to rank $m$ for different sample sizes on wavelet approximation space $\mathbb W_{4,3}\otimes \mathbb W_{4,3}$ ($P=1600$)}
\label{rastrigin_rank}
\end{figure}

Figure \ref{rastrigin_rank} shows the convergence of the approximation obtained with Algorithm \ref{alg:sparse_rank_m} using different sample sizes. We find that as the sample size increases, we get better approximations with increasing rank. {\color{black} Conversely, if only very few samples are available, then a very low rank approximation, even rank one, is able to capture the global features. The proposed method provides the simplest representation of the function with respect to the available information.}


\begin{figure}[h]
\centering
\includegraphics[width=0.55\textwidth]{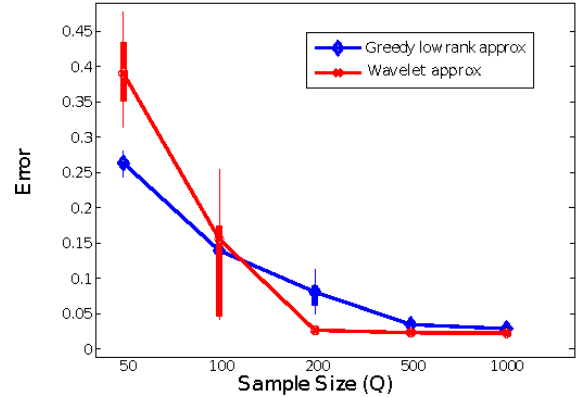}
\caption{Evolution of error $\varepsilon(u_{m_{op}},u)$ with respect to sample size $Q$. (Red line) approximation obtained with direct least-squares approximation with $\ell_1$-regularization on the full polynomial wavelet approximation space $\mathbb W_{4,3}\otimes \mathbb W_{4,3}$, (blue line) approximation obtained with Algorithm \ref{alg:sparse_rank_m} (with $\ell_1$-regularization) and with optimal rank selection.}
\label{fig:rastrigin_tensor_poly_comp}
\end{figure}

We finally analyze the robustness of Algorithm \ref{alg:sparse_rank_m} 
with respect to the sample sets. We use wavelet bases. An optimal rank approximation $u_{m_{op}}$ is selected using 3-fold cross validation as described in section \ref{sec:rankselection}. 
We compare this algorithm with a direct sparse least-squares approximation in the tensorized polynomial wavelet space (no low-rank approximation), using $\ell_1$-regularization (use of Algorithm \ref{alg:loo_cv}). Figure \ref{fig:rastrigin_tensor_poly_comp} shows the evolution of the relative error with respect to the sample size $Q$ for these two strategies.  
The vertical lines represent the scattering of the error when different sample sets are used. We observe a smaller variance of the obtained approximations when exploiting low-rank representations. This can be explained by the lower dimensionality of the representation, which is better estimated with a few number of samples. On this simple example, we see the interest of using greedy constructions of sparse low-rank representations when only a small number of samples is available, {\color{black}indeed the problem is reduced to one where elements of subsets of small dimension are to be learnt. The interest of using low-rank representations should also become clear when dealing with higher dimensional problems.}


\begin{figure}[h!]\centering

\subfigure[$Q=80$]{\includegraphics[width=8cm]{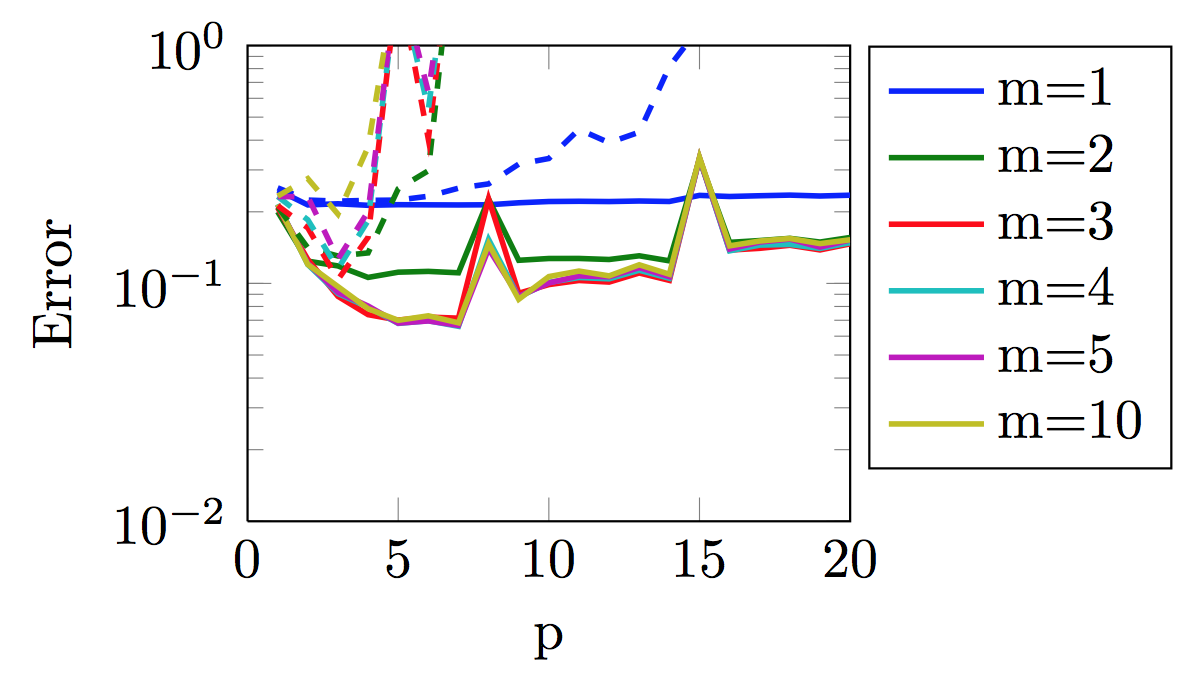}\includegraphics[width=8cm]{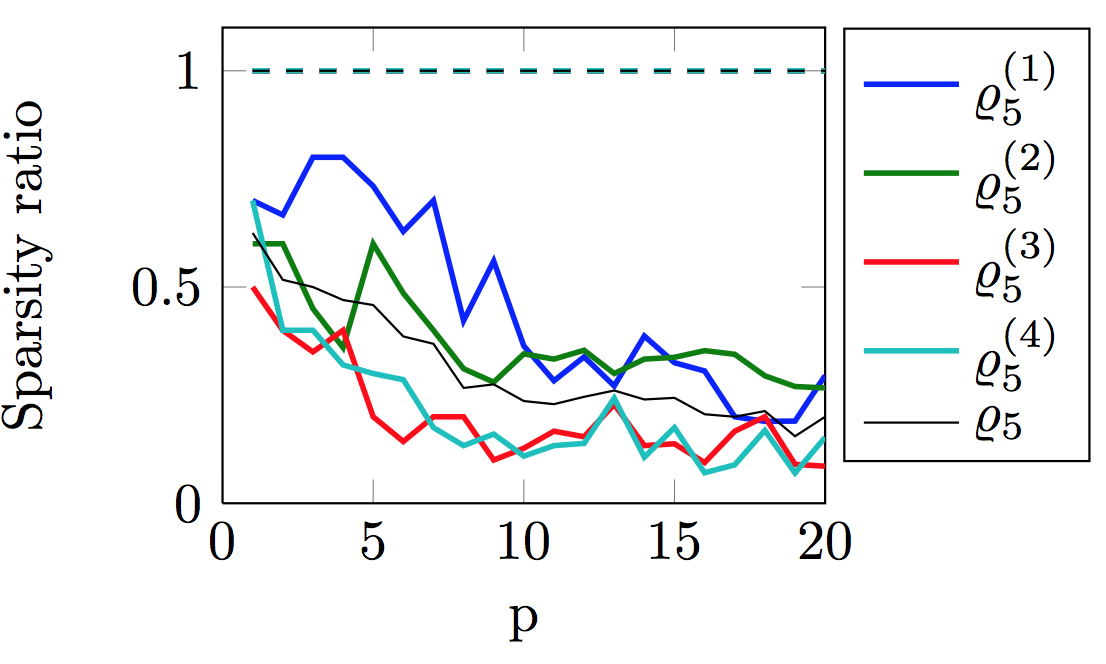}}\\
\subfigure[$Q=200$]{\includegraphics[width=8cm]{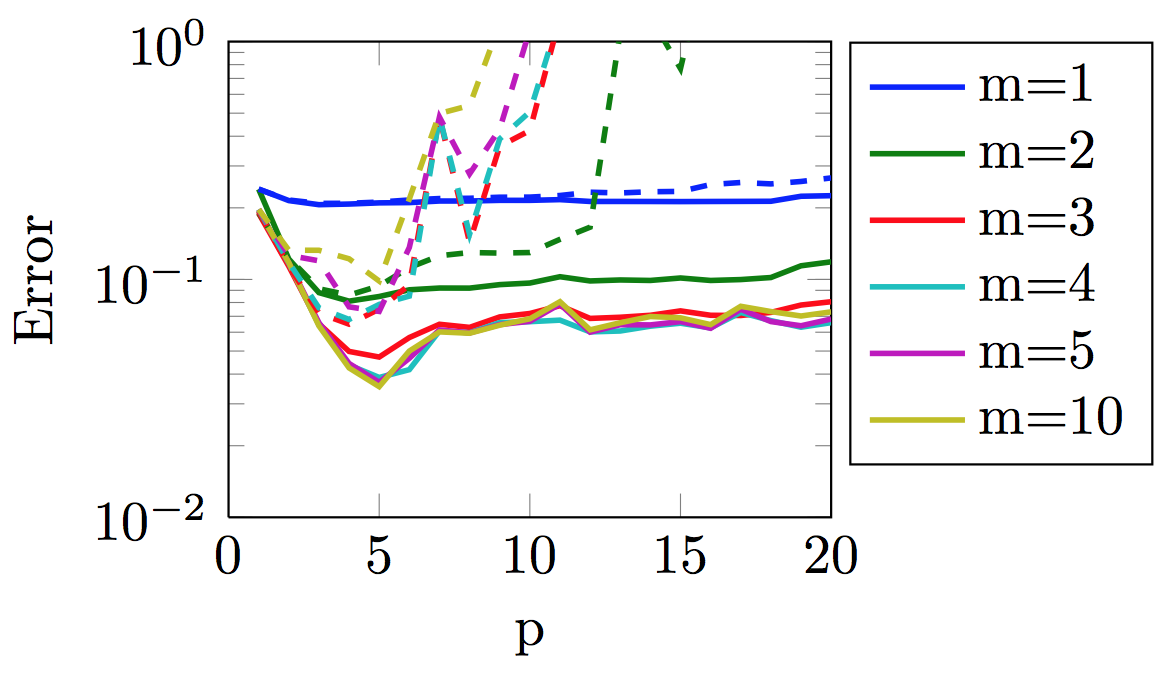}\includegraphics[width=8cm]{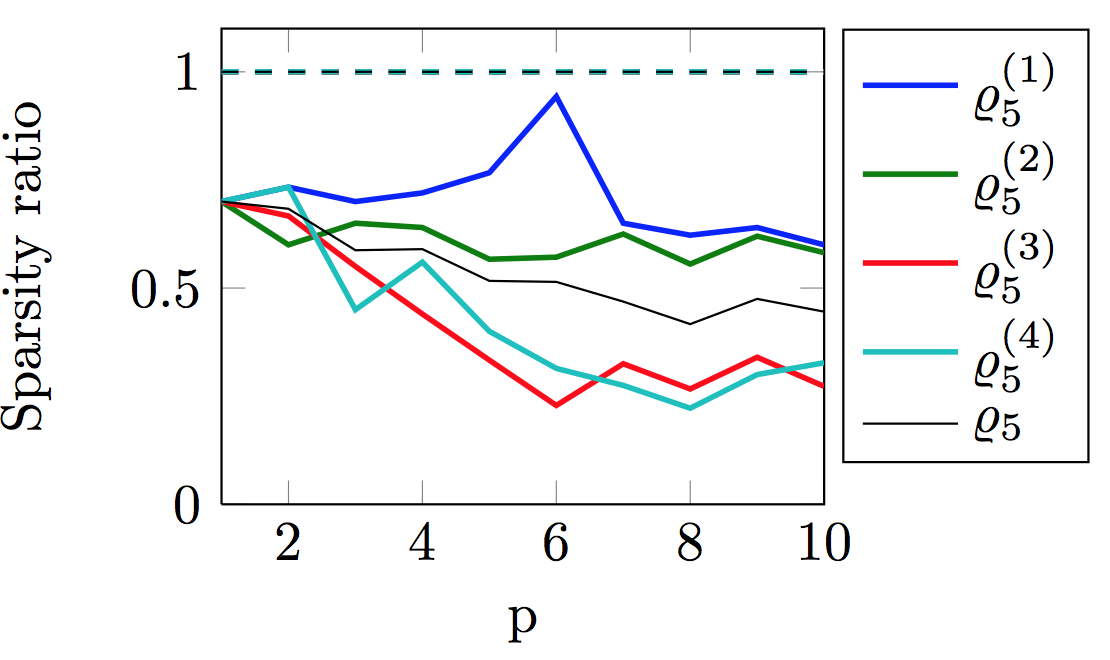}}\\
\subfigure[$Q=500$]{\includegraphics[width=8cm]{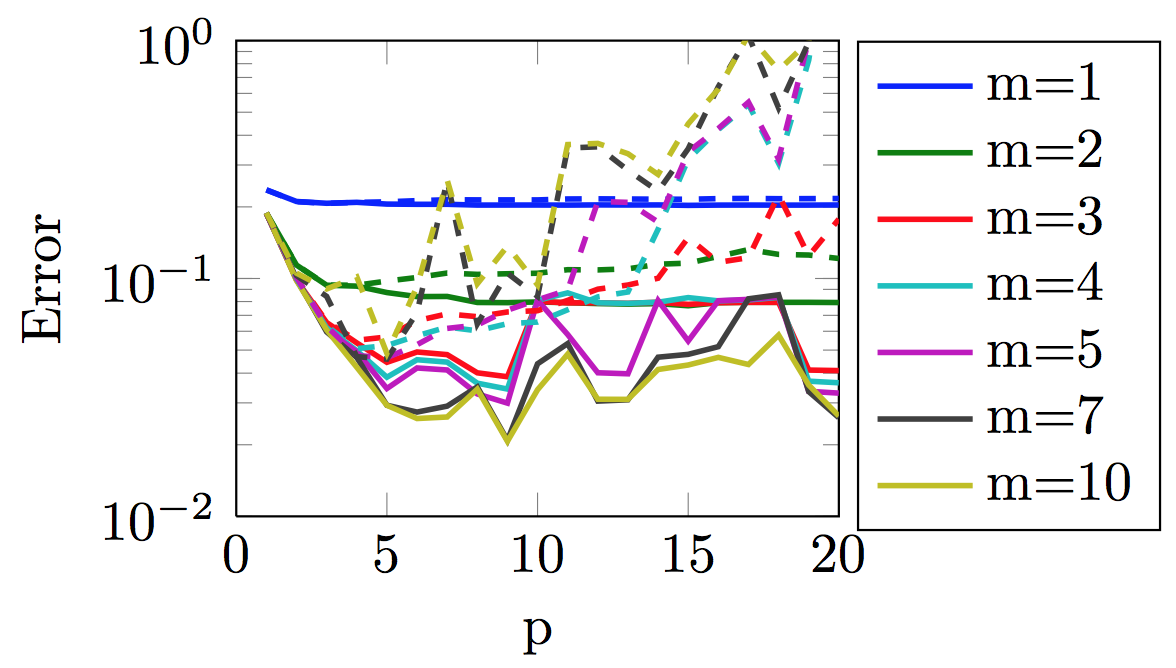}\includegraphics[width=8cm]{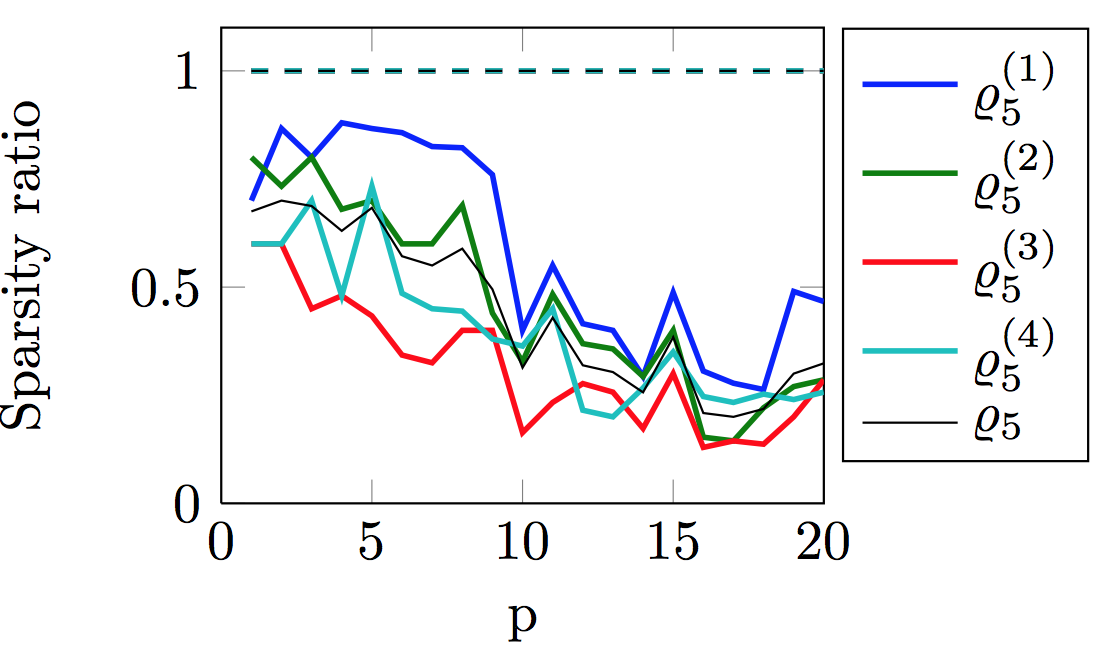}}
  \caption{Left column: evolution of approximation error $\varepsilon(I_m,I)$ with respect to polynomial degree $p$ with a fixed number of samples $Q=80,Q=200$ and $Q=500$. Polynomial approximations obtained using $\ell_1$-regularization (solid lines) and using OLS (dashed lines). Right column: total and partial sparsity ratios with respect to polynomial degree for $m=5$. }
  \label{fig:twoplates_l1ols_gpc}
\end{figure}

\begin{figure}[h!]\centering

\includegraphics[width=8cm]{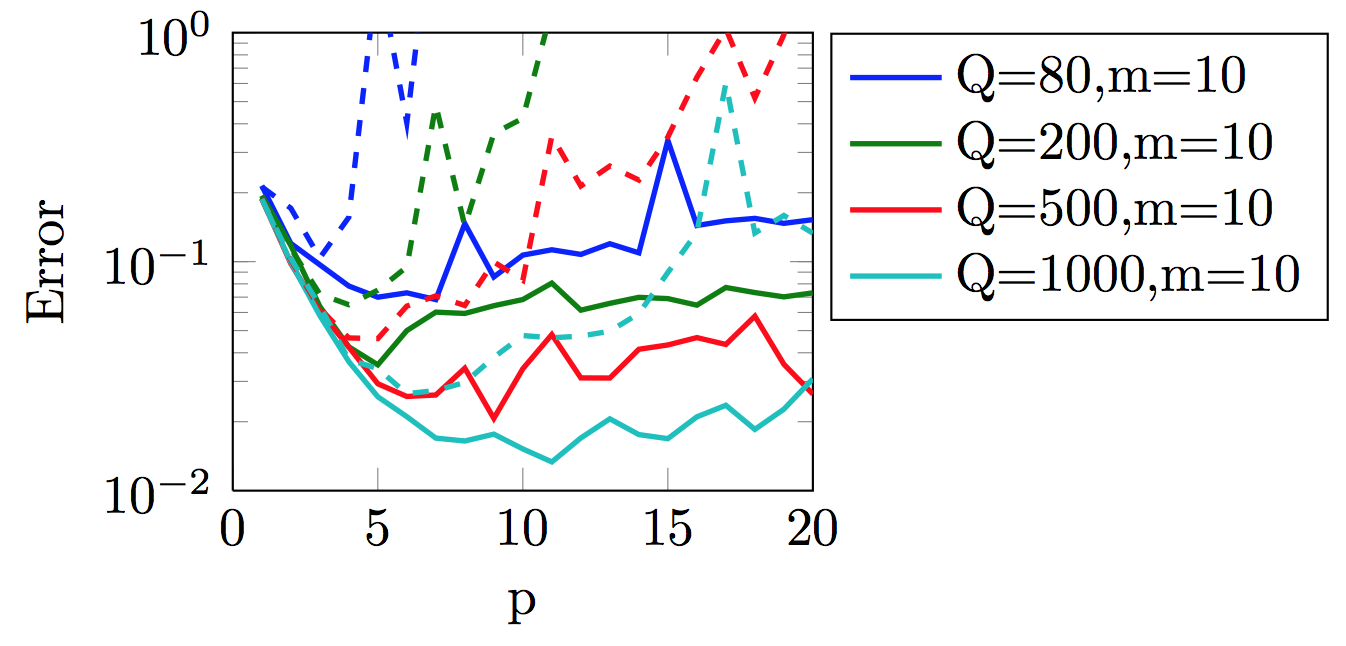}
  \caption{Evolution of approximation error $\varepsilon(I_m,I)$ with respect to polynomial degree $p$ for different sizes of sample sets (random sampling). Polynomial approximations obtained using $\ell_1$-regularization (solid lines) and using OLS (dashed lines).}
  \label{fig:twoplates_l1ols_gpc_summary}
\end{figure}

\begin{figure}

	\centering
	\subfigure[$\xi_3=0,\xi_4=0$]{\includegraphics[width=5cm]{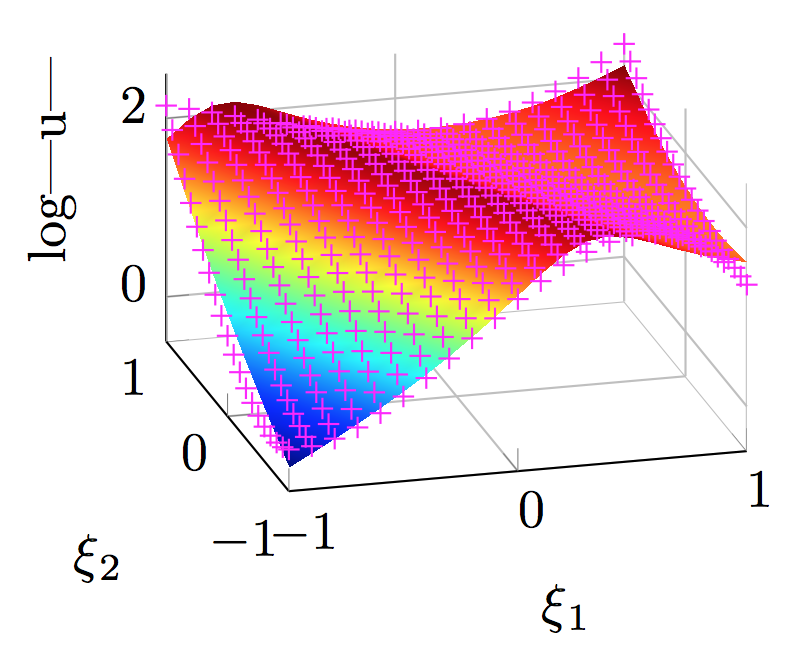}}
	\subfigure[$\xi_1=0,\xi_2=0$]{\includegraphics[width=8cm]{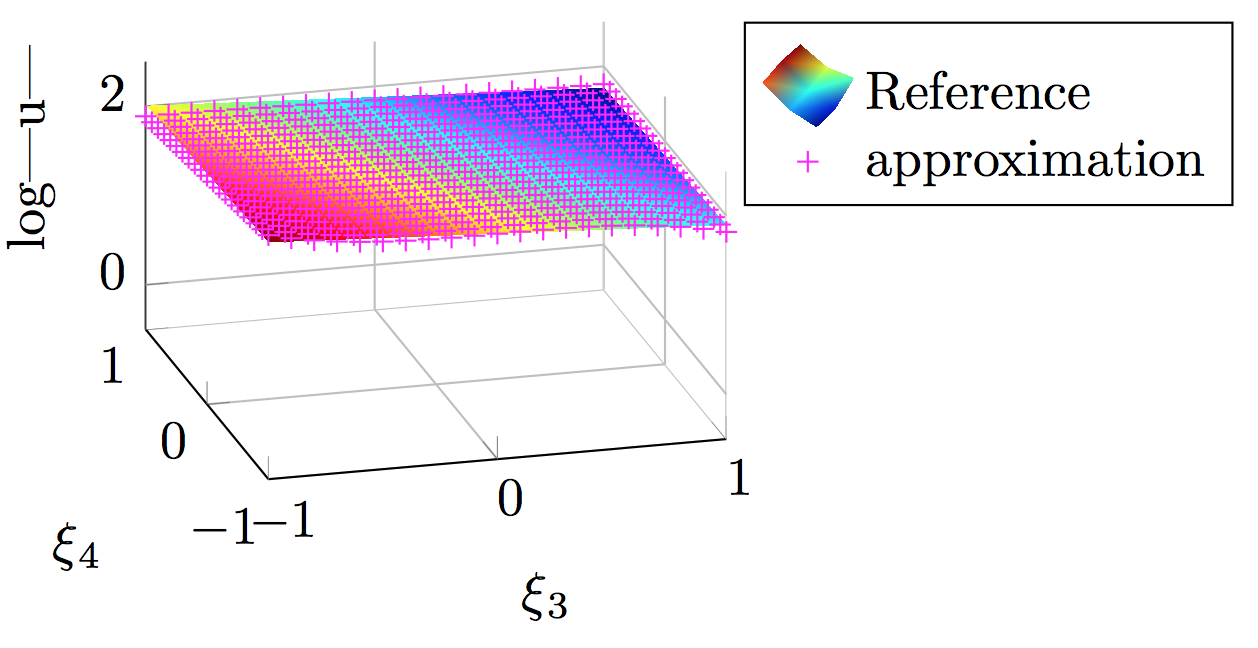}}
	\caption{Response surface : (surface) reference and (dotted) rank-4 approximation obtained with polynomial basis $p=4$, $Q=200$ using $l_1$ regularization.  }
  \label{fig:twoplates_l1_respsurf}
\end{figure}

{\color{black}
\subsection{A model problem in structural vibration analysis}\label{sec:twoplates}
\subsubsection{Model problem}

We consider a forced vibration problem of a slightly damped random linear elastic structure. The structure composed of two plates is clamped on  part $\Gamma_1$  of the boundary and submitted to a harmonic load on part  $\Gamma_2$ of the boundary as represented in figure \ref{fig:twoplates_model}(a). A finite element approximation is introduced at the spatial level using a mesh composed of $1778$ DKT
plate elements (see figure \ref{fig:twoplates_model}(b)) and leading to a discrete deterministic model with $N=5556$ degrees of freedom.
\begin{figure}[!htb]
\centering
\subfigure[]{\includegraphics[width=4cm]{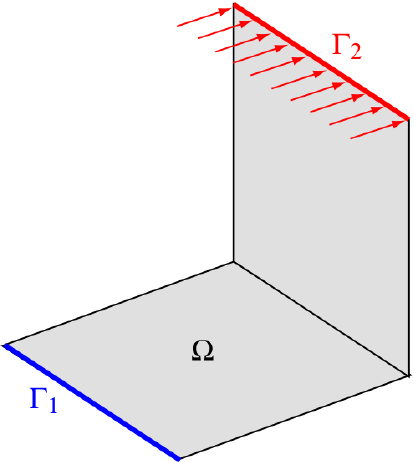}}
\qquad\subfigure[]{\includegraphics[width=4cm]{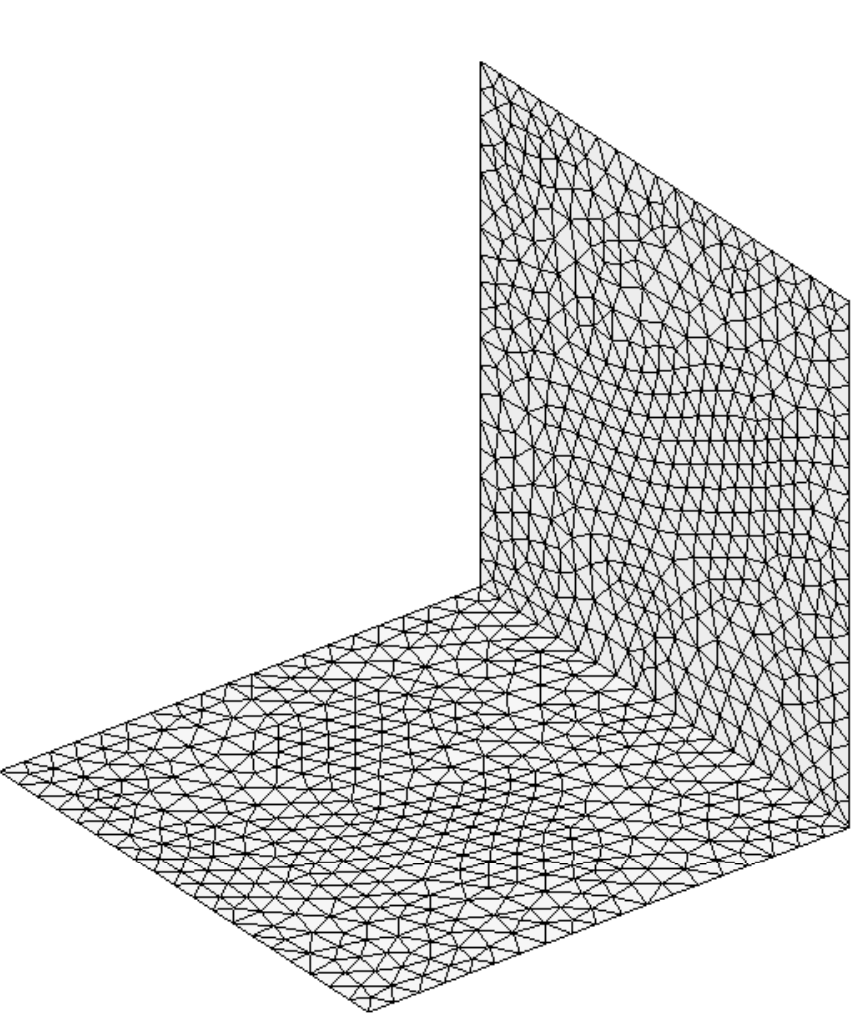}}
\caption{Elastic plate structure under
harmonic bending load. Geometry and boundary conditions (a) and
finite element mesh (b).} \label{fig:twoplates_model}
\end{figure}

The resulting discrete problem at frequency $\omega$ writes
\begin{align*}
(-\omega^2 \Mbsf+ i\omega \Cbsf+\Kbsf)\ubsf= \fbsf, \label{eq:problem}
\end{align*}
where $\ubsf\in \Cbb^N$ is the vector of coefficients of the approximation of the displacement field and $\Mbsf$, $\Kbsf=E \tilde \Kbsf$ and $\Cbsf=i\omega \eta E \tilde \Kbsf$ are the mass, stiffness and damping matrices respectively. The non-dimensional analysis considers a unitary mass density and a circular frequency $\omega=0.67 \text{ rad.s}^{-1}$. The Young modulus $E$ and the damping parameter $\eta$ are defined by 
\begin{align*}
&E=
\begin{cases} 
0.975+0.025\xi_1 \qquad\text{ on horizontal plate,} \\
0.975+0.025\xi_2 \qquad\text{ on vertical plate,} \\
\end{cases}\\
&\eta=
\begin{cases} 
0.0075+0.0025\xi_3 \quad \text{ on horizontal plate,} \\
0.0075+0.0025\xi_4 \quad \text{ on vertical plate,} \\
\end{cases}
\end{align*}
where the $\xi_k \sim U(-1,1)$, $k={1,\cdots,4}$, are
independent uniform random variables (here $d=4$). We define the quantity of
interest
\begin{equation}
I(u)(\xi)=\log\norm{u_{c}},  \notag
\end{equation}
where $u_c$ is the displacement of the top right node of the two plate structure.

\subsubsection{Impact of regularization and stochastic polynomial degree}
In this example, we illustrate that the approximation in sparse low rank tensor subsets is robust when increasing the degree of underlying polynomial approximation spaces with a fixed number of samples $Q$. 

We use Legendre polynomial basis functions with degree $p=1$ to $20$ and denote by $\Pbb_p$ the corresponding space of polynomials of maximal degree $p$ in each dimension. A rank-$m$ approximation is searched in the isotropic tensor space $\Pbb_p\otimes \cdots \otimes \Pbb_p$.

Figure \ref{fig:twoplates_l1ols_gpc}(left column) shows the error $\varepsilon(I_m,I)$ as a function of the polynomial degree $p$ for different ranks $m$  
for three different sizes of the sample set, $Q=80,200$ and $500$. 
The low rank approximation $I_m$ is computed with and without sparsity constraint, $i.e.$ we compare OLS (dashed lines) and $\ell_1$ regularization (solid lines) in correction step \ref{algo3-correction} of Algorithm \ref{alg:sparse_rank_m}.
Figure \ref{fig:twoplates_l1ols_gpc_summary} summarizes the error $\varepsilon(I_m,I)$ for different sizes of sample sets for the rank-$10$ approximation when using $\ell_1$-regularization (solid lines) and for the rank-$m$ approximation giving the best approximation when using OLS (dashed lines). 
On the one hand, we find that OLS yields a deterioration with higher polynomial order. This is consistent with the conclusions in section \ref{sec:friedman} and the quadratic rule according to which convergence is observed for $Q\geq dm(p+1)^2 $ and a deterioration is expected otherwise.
On the other hand, we see that $\ell_1$-regularization gives a more stable approximation with increasing polynomial order and also gives a more accurate best approximation than the best approximation obtained with OLS. This can be attributed to the selection of pertinent basis functions obtained by imposing sparsity constraint. 
Indeed, we clearly see in figure \ref{fig:twoplates_l1ols_gpc}(right column) that the sparsity ratio $\rho_5$ for sparse low rank approximation (solid black line) decreases with increasing polynomial degree. 
Along with the total sparsity ratio, the partial sparsity ratios $\varrho_5^{(k)}$ in each dimension $k=1$ to $4$ are plotted in figure \ref{fig:twoplates_l1ols_gpc}(right column). 
We see that $\ell_1$-regularization exploits sparsity especially in dimensions $3$ and $4$ corresponding to the damping coefficients, indeed the quantity of interest has smooth dependance on variables $\xi_3$ and $\xi_4$ whereas it has a high non linear behavior with respect to $\xi_1$ and $\xi_2$.
Figure \ref{fig:twoplates_l1_respsurf} shows the reference quantity of interest and the rank-$3$ approximation $I_3$ obtained using $\ell_1$-regularization and polynomial degree $5$ constructed from $Q=200$ samples.

This illustration also points out that a small number of model evaluations, for instance $Q=80$, does not enable to capture correctly local features of the function and a low rank approximation ($m=3$) is selected as the best approximation regarding the available information. As the number of samples increases, higher rank approximation are selected that capture the local features of the function more accurately.

}

\section{Conclusion}
A non-intrusive least-squares-based sparse low-rank tensor approximation method  has been proposed for propagation of uncertainty in
high dimensional stochastic models.
Greedy algorithms for low-rank tensor approximation have been combined with sparse least-squares approximation methods in order to obtain a robust construction of sparse low-rank tensor approximations in high dimensional approximation spaces when having only very few information on the function.  
The ability of the proposed
method to detect and exploit low-rank and sparsity  
was illustrated on three analytical models and on a  partial differential
equation with random coefficients. {\color{black}In order to  exploit at best the few samples available in practical applications in uncertainty propagation, algorithms that are able to automatically detect low-rank structures of a function (e.g. by finding an optimal tree in hierarchical tensor representations) should be developed.}
Also, the design of sampling strategies that are adapted to the construction of low-rank approximations could further improve the performances of these techniques.


\begin{thebibliography}{10}

\bibitem{BAC12}
{\sc F.~Bach, R.~Jenatton, J.~Mairal, and G.~Obozinski}, {\em Optimization with
  sparsity-inducing penalties}, Foundations and Trends in Machine Learning, 4
  (2012), pp.~1--106.

\bibitem{BEL11}
{\sc G.~Beylkin, B.~Garcke, and M.J. Mohlenkamp}, {\em Multivariate regression
  and machine learning with sums of separable functions}, Journal of
  Computational Physics, 230 (2011), pp.~2345--2367.

\bibitem{BLA11}
{\sc G.~Blatman and B.~Sudret}, {\em Adaptive sparse polynomial chaos expansion
  based least angle regression}, Journal of Computational Physics, 230 (2011),
  pp.~2345--2367.

\bibitem{CAN06b}
{\sc EJ. Candes, J.~Romberg, and T.~Tao}, {\em {Near optimal signal recovery
  from random projections: Universal encoding strategies?}}, IEEE Transactions
  on information theory, 52(12) (2006), pp.~5406--5425.

\bibitem{CAW04}
{\sc G.C. Cawley and N.L.C. Talbot}, {\em Fast exact leave-one-out
  cross-validation of sparse least-squares support vector machines}, Neural
  Networks, 17 (2004), pp.~1467--1475.

\bibitem{CHE99}
{\sc S.~S. Chen, D.~L. Donoho, and M.~A. Saunders}, {\em Atomic decomposition
  by basis pursuit}, SIAM Journal on Scientific Computing, 20 (1999),
  pp.~33--61.

\bibitem{DON06}
{\sc D.L. Donoho}, {\em {Compressed Sensing}}, IEEE Transactions on information
  theory, 52(4) (2006), pp.~1289--1306.

\bibitem{DOO09}
{\sc A.~Doostan and G.~Iaccarino}, {\em A least-squares approximation of
  partial differential equations with high-dimensional random inputs}, Journal
  of Computational Physics, 228 (2009), pp.~4332--4345.

\bibitem{DOO11}
{\sc A.~Doostan and H.~Owhadi}, {\em A non-adapted sparse approximation of pdes
  with stochastic inputs}, Journal of Computational Physics, 230 (2011),
  pp.~3015--3034.

\bibitem{DOO12}
{\sc A.~Doostan, A.~Validi, and G.~Iaccarino}, {\em Non-intrusive low-rank
  separated approximation of high-dimensional stochastic models},
  http://arxiv.org/abs/1210.1532v1,  (2012).

\bibitem{EFR04}
{\sc B.~Efron, T.~Hastie, I.~Johnstone, and R.~Tibshirani}, {\em Least angle
  regression}, The Annals of Statistics, 32 (2004), pp.~407--499.

\bibitem{Falco:2012fk}
{\sc A.~Falc{\'o} and A.~Nouy}, {\em Proper generalized decomposition for
  nonlinear convex problems in tensor banach spaces}, Numerische Mathematik,
  121 (2012), pp.~503--530.

\bibitem{GHA91}
{\sc R.~Ghanem and P.~Spanos}, {\em Stochastic finite elements: a spectral
  approach}, Springer, Berlin, 1991.

\bibitem{GRA13}
{\sc L.~Grasedyck, D.~Kressner, and C.~Tobler}, {\em A literature survey of
  low-rank tensor approximation techniques}, arXiv:1302.7121,  (2013).

\bibitem{HAC12}
{\sc W.~Hackbusch}, {\em Tensor Spaces and Numerical Tensor Calculus}, vol.~42
  of Series in Computational Mathematics, Springer, 2012.

\bibitem{KHO12}
{\sc B.N. Khoromskij}, {\em Tensor structured numerical methods in scientific
  computing: Survey on recent advances}, Chemometrics and Intelligent
  laboratory Systems, 110 (2012), pp.~1--19.

\bibitem{KHO10}
{\sc B.N. Khoromskij and C.~Schwab}, {\em Tensor-structured galerkin
  approximation of parametric and stochastic elliptic pdes}, Tech. Report
  Research Report No. 2010-04, ETH, 2010.

\bibitem{KOL09}
{\sc T.~G. Kolda and B.~W. Bader}, {\em Tensor decompositions and
  applications}, SIAM Review, 51 (2009), pp.~455--500.

\bibitem{LEM04}
{\sc O.P. {Le~Ma\^{i}tre}, O.~M. Knio, H.~N. Najm, and R.~G. Ghanem}, {\em
  Uncertainty propagation using {W}iener-{H}aar expansions}, Journal of
  Computational Physics, 197 (2004), pp.~28--57.

\bibitem{LEM10}
{\sc O.~P. {Le Ma\^{i}tre} and O.~M. Knio}, {\em Spectral Methods for
  Uncertainty Quantification With Applications to Computational Fluid
  Dynamics}, Scientific Computation, Springer, 2010.

\bibitem{MAI10}
{\sc J.~Mairal, F.~Bach, J.~Ponce, and G.~Sapiro}, {\em Online learning for
  matrix factorization and sparse coding}, Journal of Machine Learning
  Research, 11 (2010).

\bibitem{MAT08}
{\sc H.~G. Matthies}, {\em Stochastic finite elements: Computational approaches
  to stochastic partial differential equations}, {Zamm-Zeitschrift Fur
  Angewandte Mathematik Und Mechanik}, {88} ({2008}), pp.~{849--873}.

\bibitem{MAT11b}
{\sc H.~G. Matthies and E~Zander}, {\em Solving stochastic systems with
  low-rank tensor compression}, Linear Algebra and its Applications, 436
  (2012), pp.~3819--3838.

\bibitem{NOB11}
{\sc G.~Migliorati, F.~Nobile, E.~von Schwerin, and R.~Tempone}, {\em Analysis
  of the discrete {L}2 projection on polynomial spaces with random
  evaluations}, Tech. Report~46, MATHICSE, 2011.

\bibitem{NOU07}
{\sc A.~Nouy}, {\em A generalized spectral decomposition technique to solve a
  class of linear stochastic partial differential equations}, Computer Methods
  in Applied Mechanics and Engineering, 196 (2007), pp.~4521--4537.

\bibitem{NOU09d}
\leavevmode\vrule height 2pt depth -1.6pt width 23pt, {\em Recent developments
  in spectral stochastic methods for the numerical solution of stochastic
  partial differential equations}, Archives of Computational Methods in
  Engineering, 16 (2009), pp.~251--285.

\bibitem{NOU10}
\leavevmode\vrule height 2pt depth -1.6pt width 23pt, {\em Proper generalized
  decompositions and separated representations for the numerical solution of
  high dimensional stochastic problems}, Archives of Computational Methods in
  Engineering, 17 (2010), pp.~403--434.

\bibitem{RAI12}
{\sc P.~Rai, M.~Chevreuil, A.~Nouy, and J.~Sen Gupta}, {\em A regression based
  non-intrusive method using separated representation for uncertainty
  quantification}, in ASME 2012, 11th Biennial Conference on Engineering
  Systems Design and Analysis (ESDA 2012), Nantes, France, July 2-4 2012.

\bibitem{TIB96}
{\sc R.~Tibshirani}, {\em Regression shrinkage and selection via the lasso},
  Journal of the Royal Statistical Society Series B, 58 (1996), pp.~267--288.

\bibitem{XIU02}
{\sc D.~Xiu and G.~E. Karniadakis}, {\em The {W}iener-{A}skey polynomial chaos
  for stochastic differential equations}, SIAM J. Sci. Comput., 24 (2002),
  pp.~619--644.

\end{thebibliography}

\end{document}